\documentclass[12pt]{report}
\usepackage{aums}       % For Master's papers
\usepackage{ulem}       % underlining on style-page; see \normalem below
\usepackage{url}
\usepackage{tikz}
\usepackage{pgf}
\usepackage[intoc]{nomencl}
   	       
\makenomenclature 
\usepackage{enumerate}
\usepackage{amssymb, amsmath}
\usepackage{amsthm}
\usepackage{cleveref}

\newcounter{dummy} \numberwithin{dummy}{chapter}
\newtheorem{theorem}[dummy]{Theorem}
\newtheorem{prop}[dummy]{Proposition}
\newtheorem{lemma}[dummy]{Lemma}
\newtheorem{corollary}[dummy]{Corollary}

\title{Matrix Algebras over Strongly Non-Singular Rings }
\author{Bradley McQuaig} 
\date{May 4, 2014} %date of graduation
\copyrightyear{2014} %copyright year

\keywords{strongly non-singular, torsion-free, Baer-ring, semi-hereditary}

\adviser{Ulrich Albrecht}

\professor{Ulrich Albrecht, Professor of Mathematics}
\professor{H. Pat Goeters, Professor of Mathematics}
\professor{Georg Hetzer, Professor of Mathematics}

\begin{document}

\begin{romanpages}      % roman-numbered pages 

\TitlePage 

\begin{abstract}
We consider some existing results regarding rings for which the classes of torsion-free and non-singular right modules coincide.
Here, a right $R$-module $M$ is {\it non-singular} if $xI$ is nonzero for every nonzero $x \in M$ and every essential right ideal $I$ of $R$, and a right $R$-module $M$ is {\it torsion-free} if $Tor_{1}^{R}(M,R \slash Rr)=0$ for every $r \in R$.
In particular, we consider a ring $R$ for which the classes of torsion-free and non-singular right $S$-modules coincide for every ring $S$ Morita-equivalent to $R$.
We make use of these results, as well as the existence of a Morita-equivalence between a ring $R$ and the $n \times n$ matrix ring $Mat_{n}(R)$, to characterize rings whose $n \times n$ matrix ring is a Baer-ring. 
A ring is {\it Baer} if every right (or left) annihilator is generated by an idempotent.
Semi-hereditary, strongly non-singular, and Utumi rings will play an important role, and we explore these concepts and relevant results as well.

\end{abstract}

%\begin{acknowledgments}
%Put text of the acknowledgments here.
%\end{acknowledgments}

\tableofcontents
%\listoffigures
%\listoftables

\printnomenclature[0.5in] %used for the List of Abbreviations
\end{romanpages}        % All done with roman-numbered pages

%\normalem       % Make italics the default for \em

\chapter{Introduction}  % Use \\ for long titles 

% This is a sample document for the Auburn \nomenclature{Auburn}{Auburn University} \LaTeX{} style-files known
%as {\tt aums} (for Master's papers) and {\tt auphd} (for Ph.D.'s).
%The appendix contains some of the history of this project, including
%contact information for the authors. 

%Site administrators should upgrade to \LaTeX2e; however, the style files
%should work with the older \LaTeX.  The style files should be available
%on mallard.  The current release is available by anonymous ftp to
%ftp.dms.auburn.edu in the directory aums (on-campus computers may
%also retrieve these from \url{http://www.dms.auburn.edu/manuals}).
%Most users will need either Lamport's book  or Hahn's book 

%If you do not need the List of Abbreviations\nomenclature{LoA}{List of Abbreviations}, comment the nomencl package and associated nomenclature commands. 

In this thesis, we consider the relationship between a ring $R$ and $Mat_{n}(R)$, the $n \times n$ matrix ring over $R$.
In particular, we investigate necessary and sufficient conditions placed on $R$ so that $Mat_{n}(R)$ is a Baer-ring. 
A ring is a {\it Baer-ring} if every right (or left) annihilator ideal is generated by an idempotent.
In determining these conditions, we make use of the existence of a Morita-equivalence between $R$ and $Mat_{n}(R)$ (\Cref{mat morita}), as well as the fact that $Mat_{n}(R)$ is isomorphic to the endomorphism ring of any free right $R$-module with basis $\{x_{i}\}_{i=1}^{n}$(\Cref{End iso Mat}).
Here, two rings are {\it Morita-equivalent} if their module categories are equivalent, and the {\it endomorphism ring} $End_{R}(M)$ of a right $R$-module $M$ is the set of all $R$-homomorphisms $f:M \rightarrow M$, which is a ring under pointwise addition and composition of functions.

The concepts of torsion-freeness and non-singularity of modules will also come into play.
In particular, we consider rings for which the classes of torsion-free and  non-singular right $S$-modules coincide for every ring $S$ Morita-equivalent to $R$.
Albrecht, Dauns, and Fuchs investigate such rings in \cite{Albrecht}.
A module $M$ over a ring $R$ is {\it torsion-free in the classical sense} if $xr \neq 0$ for every nonzero $x \in M$ and every regular $r \in R$, where $r \in R$ is {\it regular} if it is not a left or right zero-divisor.
For commutative rings, this is a useful way to define such modules, especially for integral domains since regular elements are precisely the nonzero elements.
In the case $R$ is non-commutative, then the set $M_{t}=\{x \in M$ $|$ $ann_{r}(x)$ contains some regular element of R$\}$, which is usually referred to as the {\it torsion-submodule} in the commutative setting, is not necessarily a submodule of $M$.
There are other ways in which torsion-freeness can be defined in the non-commutative setting.
In \cite{Hattori}, Hattori calls a right $R$-module $M$ {\it torsion-free} if $Tor_{1}^{R}(M,R/Rr)=0$ for every $r \in R$.  
This is based on homological properties of modules and coincides with the classical definition in the case $R$ is commutative. 
In \cite{Goodearl}, Goodearl defines the singular submodule and non-singularity of modules in the general non-commutative setting, which is closely related to the concept of torsion submodules and torsion-freeness.
We look at relevant background information on torsion-freeness and non-singularity in Chapters 4 and 5.

Albrecht, Dauns, and Fuchs found that $S$ is right strongly non-singular and the classes of torsion-free and non-singular $S$-modules coincide for every ring $S$ Morita equivalent to a ring $R$ if and only if $R$ is right strongly non-singular, right semi-hereditary, and does not contain an infinite set of orthogonal idempotents \cite[Theorem 5.1]{Albrecht}.
A ring is {\it right strongly non-singular} if its maximal right ring of quotients is a perfect left localization.
These rings will be explored in Section 5.2, and semi-hereditary rings will be defined and explored in Chapter 2.
We make use of this theorem and take it a step further to show that $Mat_{n}(R)$ is a right and left Utumi Baer-ring if and only if the classes of torsion-free and non-singular $S$-modules coincide for every ring $S$ Morita equivalent to a ring $R$.
Note that we remove the condition that every Morita-equivalent ring $S$ need be strongly non-singular.
Instead, we assume that our ring $R$ is right Utumi, and from this we also get that $Mat_{n}(R)$ is both right and left Utumi.
We define Utumi rings in Section 5.3.

Unless noted otherwise, commutativity of a ring is not assumed, but all rings are assumed to have a multiplicative identity.

\chapter{Semi-hereditary Rings and p.p.-rings}

We begin by looking at projective modules.
A right R-module P is {\it projective} if given right R-modules A and B, an epimorphism $\pi : A \rightarrow B$, and a homomorphism $\varphi : P \rightarrow$ B, then there exists a homomorphism $\psi$ : P $\rightarrow$ A such that $\pi\psi$ = $\varphi$.
In particular, every free right $R$-module is projective \cite[Theorem 3.1]{Rotman}.
We make use of the following well-known characterization of projective modules:

\begin{theorem}\cite{Rotman}\label{projective}
Let $R$ be a ring.  The following are equivalent for a right $R$-module $P$:
\begin{enumerate}[(a)]
   \item $P$ is projective
   \item $P$ is isomorphic to a direct summand of a free right $R$-module.  In other words, there is a free right $R$-module $F = Q \bigoplus N$, where $N$ is a right $R$-module and $Q \cong P$. 
   \item For any right $R$-module $M$ and epimorphism $\varphi:M \rightarrow P$, $M = \ker{(\varphi)} \bigoplus N$.
  \end{enumerate}
\end{theorem}

Let $Mod_{R}$ be the category of all right R-modules for a ring R.  
A {\it complex} in $Mod_{R}$ is a sequence of right R-modules and R-homomorphisms in $Mod_{R}$, 
\newline \centerline{$... \rightarrow A_{k+1} \xrightarrow{\alpha_{k+1}} A_{k} \xrightarrow{\alpha_{k}} A_{k-1} \rightarrow ...$}
\newline such that $\alpha_{k+1}\alpha_{k} = 0$ for every $k \in \mathbb{Z}$.
Observe $\alpha_{k+1}\alpha_{k} = 0$ implies that $im(\alpha_{k+1}) \subseteq \ker{(\alpha_{k})}$.
The sequence is called {\it exact} if $im (\alpha_{k+1}) = \ker{(\alpha_{k})}$ for every $k \in \mathbb{Z}$.
An exact sequence $0 \rightarrow A \xrightarrow{\alpha} B \xrightarrow{\beta} C \rightarrow 0$ of right $R$-modules is referred to as a {\it short exact sequence}.
Such an exact sequence is said to {\it split} if there exists an $R$-homomorphism $\gamma: C \rightarrow B$ such that $\beta\gamma = 1_{C}$, where $1_{C}$ is the identity map on $C$.

\begin{lemma}\cite{Rotman}\label{Rot 2.28}
Let $0 \rightarrow A \xrightarrow{\alpha} B \xrightarrow{\beta} C \rightarrow 0$ be a sequence of right $R$-modules.
If this sequence is split exact, then $B \cong A \bigoplus C$.
\begin{proof}
If the exact sequence $0 \rightarrow A \xrightarrow{\alpha} B \xrightarrow{\beta} C \rightarrow 0$ of right $R$-modules splits, then there exists an $R$-homomorphism $\gamma:C \rightarrow B$ such that $\beta\gamma \cong 1_{C}$.
Observe that since $\alpha$ is a monomorphism, $im(\alpha) \cong A$.
Moreover, if $x \in \ker(\gamma)$, then $\gamma(x)=0$.
However, $\beta(0)=\beta\gamma(x)=x$ since $\beta\gamma=1_{C}$.
Thus, $x = 0$ and $\gamma$ is also a monomorphism.
Hence, $im(\beta) \cong C$.
Therefore, to show that $B \cong A \bigoplus C$, it suffices to show that $B \cong im(\alpha) \bigoplus im(\gamma)$.

Let $b \in B$.
Then $\beta(b) \in C$ and $\gamma\beta(b) \in im(\gamma)$.
Furthermore, $b - \gamma\beta(b) \in \ker(\beta)=im(\alpha)$ since $\beta(b-\gamma\beta(b))=\beta(b)-\beta\gamma\beta(b)=\beta(b)-\beta(b)=0$.
Hence, $b=[b-\gamma\beta(b)]+\gamma\beta(b) \in im(\alpha)+im(\gamma)$.
Suppose, $x \in im(\alpha) \cap im(\gamma)$.
Then, there exists some $a \in A$ such that $\alpha(a)=x$, and there exists some $c \in C$ such that $\gamma(c)=x$.
Now, $\alpha(a) \in im(\alpha)=\ker(\beta)$, which implies $\beta(x)=\beta\alpha(a)=0$.
However, it is also the case that $\beta(x)=\beta\gamma(c)=c$.
Hence, $c=0$ and it follows that $x=\gamma(c)=\gamma(0)=0$.
Thus, $im(\alpha) \cap im(\gamma)=0$.
Therefore, $B \cong im(\alpha) \bigoplus im(\gamma) \cong A \bigoplus C$.
\end{proof}
\end{lemma}

\begin{prop}\cite{Rotman}\label{Rotman 3.2}
The following are equivalent for a right $R$-module $P$:
\begin{enumerate}[(a)]
\item $P$ is projective.
\item The sequence $0 \rightarrow Hom_{R}(P,A) \xrightarrow{Hom_{R}(P, \varphi)} Hom_{R}(P,B) \xrightarrow{Hom_{R}(P, \psi)} Hom_{R}(P,C) \rightarrow 0$ is exact whenever $0 \rightarrow A \xrightarrow{\varphi} B \xrightarrow{\psi} C \rightarrow 0$ is a an exact sequence of right $R$-modules.
\end{enumerate}
\begin{proof}
($a$) $\Rightarrow$ ($b$): Suppose $P$ is projective.
Observe that the functor $Hom_{R}(P,\underline{\hspace{.3cm}})$ is left exact \cite[Theorem 2.38]{Rotman}.
Thus, if $0 \rightarrow A \xrightarrow{\varphi} B \xrightarrow{\psi} C \rightarrow 0$ is exact, then $$0 \rightarrow Hom_{R}(P,A) \xrightarrow{Hom_{R}(P, \varphi)} Hom_{R}(P,B) \xrightarrow{Hom_{R}(P, \psi)} Hom_{R}(P,C)$$ is exact.
Therefore, it remains to be shown that $Hom_{R}(P, \psi)$ is an epimorphism.
Let $\alpha \in Hom_{R}(P,C)$.
Since $P$ is projective, there exists a homomorphism $\beta:P \rightarrow B$ such that $\alpha=\psi\beta$.
Hence, $Hom_{R}(P,\psi)(\beta)=\psi\beta=\alpha$.
Therefore, $Hom_{R}(P, \psi)$ is an epimorphism.

($b$) $\Rightarrow$ ($a$): Let $P$ be a right $R$-module and assume exactness of $0 \rightarrow A \xrightarrow{\varphi} B \xrightarrow{\psi} C \rightarrow 0$ implies exactness of $0 \rightarrow Hom_{R}(P,A) \xrightarrow{Hom_{R}(P, \varphi)} Hom_{R}(P,B) \xrightarrow{Hom_{R}(P, \psi)} Hom_{R}(P,C) \rightarrow 0$.
This implies $Hom_{R}(P, \psi)$ is an epimorphism.
Thus, if $\alpha \in Hom_{R}(P,\psi)$, then there exists some $\beta \in Hom_{R}(P,B)$ such that $Hom_{R}(P,\psi)(\beta)=\psi\beta=\alpha$.
That is, given an epimorphism $\psi:B \rightarrow C$ and a homomorphism $\alpha:P \rightarrow C$, there exists a homomorphism $\beta:P \rightarrow B$ such that $\alpha=\psi\beta$.
Therefore, $P$ is projective.
\end{proof}
\end{prop}

A ring $R$ is a {\it right p.p.-ring} if every principal right ideal is projective as a right $R$-module.
A ring $R$ is {\it right semi-hereditary} if every finitely generated right ideal is projective as a right R-module.
For a right $R$-module $M$ and any subset $S \subseteq M$, define the {\it right annihilator of $S$} in $R$ as $ann_{r}(S)$ = $\{r$ $\in$ R $|$ $xr = 0$ for every $x \in S\}$.  
The right annihilator of $S$ is a right ideal of $R$.
Similarly, the {\it left annihilator of $S$} in $R$ can be defined for a left $R$-module $M$ as $ann_{l}(S)$ = $\{r$ $\in$ R $|$ $rx = 0$ for every $x \in S\}$.
The left annihilator of $S$ is a left ideal of $R$.  
The following proposition shows that right p.p.-rings can be defined in terms of annihilators of elements and idempotents, where an idempotent is an element $e \in R$ such that $e^2=e$.

\begin{prop}\label{pp idempt}A ring R is a right p.p.-ring if and only if for every $x \in R$ there exists some idempotent e $\in$ R such that $ann_{r}(x)=eR$.
\begin{proof}
For $x \in R$, consider the function $f_{x}:R \rightarrow xR$ given by $r \mapsto xr$.
This is a well-defined epimorphism.
Then $R$ is a right p.p.-ring if and only if the principal right ideal $xR$ is projective for for every $x \in R$ if and only if $\ker{(f_{x})}$ is a direct summand of $R$ for every $x \in R$.
Observe that for each $x \in R$, $\ker{(f_{x})}=ann_{r}(x)$.
Hence, $R$ is a right p.p.-ring if and only if $ann_{r}(x)$ is a direct summand of $R$.
Note that every direct summand of $R$ is generated by an idempotent since $R \cong eR \bigoplus (1-e)R$ for any idempotent $e \in R$.
Thus, as a direct summand, $ann_{r}(x) = eR$ for some idempotent $e \in R$.
Therefore, $R$ is a right p.p.-ring if and only if for every $x \in R$ there is some idempotent $e \in R$ such that $ann_{r}(x)=eR$.
\end{proof}
\end{prop}

Let $Mat_{n}(R)$ denote the set of all $n \times n$ matrices with entries in $R$.
Under standard matrix addition and multiplication, $Mat_{n}(R)$ is a ring.
A useful characterization of semi-hereditary rings is that such rings are precisely those for which $Mat_{n}(R)$ is a right p.p.-ring for every $0 < n < \omega$.
To show this, the following two lemmas will be needed:

\begin{lemma}\cite{Rotman}\label{Rot 4.30}
A ring $R$ is right semi-hereditary if and only if every finitely generated submodule $U$ of a projective right $R$-module $P$ is projective.
\begin{proof}Suppose $R$ is right semi-hereditary and let $U$ be a submodule of a projective right $R$-module $P$.
By \Cref{projective}, $P \bigoplus N$ is free for some right R-module $N$.
Hence, $P$ is a submodule of a free module, and it follows that any submodule of $P$ is also a submodule of a free module.  
Thus, without loss of generality, it can be assumed that $P$ is a free right R-module.  
Moreover, since $U$ is finitely generated, it can be assumed that $P$ is finitely generated with basis $X = \{x_{1},x_{2},...,x_{n}\}$ for some $0<n<\omega$.
%
%Let $X = \{x_{i}\}_{i \in I}$ be a basis for $P$, and let $U = a_{1}R + ... + a_{n}R$ be a finitely generated submodule of $P$.
%Observe that for each $j = 1, ..., n$, $a_{j} = x_{j_{1}}r_{j_{1}}+ ... + x_{j_{k}}r_{j_{k}}$, where $j_{l} \in I$ and $r_{j_{l}} \in R$ for $l = 1, ..., k$.
%Hence, $supp(a_{j})= \{x_{i} \in X$ $|$ $x_{i}r_{i} \neq 0\}$ is finite for each $j$ and the right $R$-module $Y$ generated by $\cup_{j=1}^n supp(a_{j})$ is finitely generated.
%Furthermore, since $\cup_{j=1}^n supp(a_{j}) \subseteq X$ and $X$ is linearly independent, $\cup_{j=1}^n supp(a_{j})$ must be linearly independent as well.
%Hence, $Y$ is a finitely generated free right $R$-module.
%Clearly $U \subseteq Y$, and thus $U$ is a submodule of a finitely generated free right $R$-module.
%Therefore, it can be assumed that $P$ is finitely generated with basis $X = \{x_{1},x_{2},...,x_{n}\}$ for $0<n<\omega$.

Inductively, it will be shown that $U$ is a finite direct sum of finitely generated right ideals.
If $n=1$, then $P=x_{1}R \cong R$.
Since submodules of the right $R$-module $R$ are right ideals, $U$ is a finitely generated right ideal.
Suppose $n>1$ and assume $U$ is a finite direct sum of finitely generated right ideals for $k<n$.
Let $V = U \cap (x_{1}R+x_{2}R+...+x_{n-1}R)$.
Then, $V$ is a finitely generated submodule of a free right $R$-module with basis $\{x_{1},x_{2},...,x_{n-1}\}$.
By assumption, $V$ is a finite direct sum of finitely generated right ideals. 
Note that if $u \in U$, then $u=v+x_{n}r$ with $v \in V$ and $r \in R$.
This expression for $u$ is unique since $X$ is a linearly independent spanning set.
Thus, the map $\varphi: U \rightarrow R$ defined by $\varphi(u) = \varphi(v+x_{n}r) = r$ is a well-defined homomorphism.

Now, $im (\varphi)$ is a finitely generated right ideal of $R$ since it is the epimorphic image of the finitely generated right $R$-module $U$.
Hence, $im (\varphi)$ is projective since $R$ is right semi-hereditary. 
Consider the short exact sequence $0 \rightarrow K \xrightarrow{\iota} U \xrightarrow{\varphi} im (\varphi) \rightarrow 0$, where $K=\ker{\varphi}$ and $\iota$ is the inclusion map.
This sequence splits since $im (\varphi)$ is projective, and thus $U \cong K \bigoplus im (\varphi)$ by \Cref{Rot 2.28}.
Hence, $U$ is a finite direct sum of finitely generated right ideals since both $K$ and $im (\varphi)$ are finitely generated right ideals. 
Since $R$ is right semi-hereditary, each of these right ideals is projective.
Therefore, $U$ is projective as the direct sum of projective right ideals.

Conversely, suppose that if $P$ is a projective right $R$-module, then every finitely generated submodule $U$ of $P$ is projective.
Let $I$ be a finitely generated right ideal of $R$.
Note that $R$ is a free right $R$-module and thus projective.
Hence, $I$ is a finitely generated submodule of $R$, and by assumption $I$ is projective.
Therefore, $R$ is right semi-hereditary.
\end{proof}
\end{lemma}

\begin{lemma}\label{End iso Mat}
Let $R$ be a ring, and $F$ a finitely generated free right $R$-module with basis $\{x_{i}\}_{i=1}^{n}$ for  $0<n<\omega$.
Then, $Mat_{n}(R) \cong End_{R}(F)$. 
\begin{proof}
Let $S=End_{R}(F)$ and take $f \in S$.
Then, $f(x_{k}) \in F$ for each $k=1,2,...,n$.
Hence, $f(x_{k})$ is of the form $\displaystyle \sum_{i=1}^{n} x_{i}a_{ik}$, where $a_{ik} \in R$ for every $i$ and every $k$.
Let $A=\{a_{ik}\}$ be the $n \times n$ matrix whose $i$-$k$th entry is $a_{ik}$, and let $\varphi:S \rightarrow Mat_{n}(R)$ be defined by $f \mapsto A$.
If $f,g \in S$ are such that $f=g$, then $f(x_{k})=g(x_{k})$ for every $k=1,2,...,n$.
Hence, $\varphi$ is well-defined.
Furthermore, if $f(x_{k})=\displaystyle \sum_{i=1}^{n} x_{i}a_{ik}$ and $g(x_{k})=\displaystyle \sum_{i=1}^{n} x_{i}b_{ik}$ for $k=1,2,...,n$, then $(f+g)(x_{k})=f(x_{k})+g(x_{k})=\displaystyle \sum_{i=1}^{n} x_{i}(a_{ik}+b_{ik})$.
Thus, if $A=\{a_{ik}\}$ and $B=\{b_{ik}\}$ are the $n \times n$ matrices with entries determined by $f$ and $g$ respectively, then $A+B=\{a_{ik}+b_{ik}\}$ is the $n \times n$ matrix with entries determined by $f+g$.
Hence, $\varphi(f+g)=A+B=\varphi f +\varphi g$.

To see that $\varphi$ is a ring homomorphism, it remains to be seen that $\varphi(fg)=\varphi(f)\varphi(g)=AB$.  
In other words, it needs to be shown that the entries of the matrix $AB$ are determined by $fg(x_{j})$ for $j=1,2,...,n$.
Observe that if $A=\{a_{ik}\}$ and $B=\{b_{ik}\}$ are $n \times n$ matrices, then under standard matrix multiplication $AB$ is the $n \times n$ matrix whose $i$-$j$th entry is $\displaystyle \sum_{k=1}^{n} a_{ik}b_{kj}$.
This is indeed the matrix determined by the endomorphism $fg$ since the following holds: \\$fg(x_{j})=f(\displaystyle \sum_{k=1}^{n} x_{k}b_{kj})=\displaystyle \sum_{k=1}^{n} f(x_{k})b_{kj}=\displaystyle \sum_{k=1}^{n}\displaystyle \sum_{i=1}^{n}x_{i}a_{ik}b_{kj}=\displaystyle \sum_{i=1}^{n}x_{i}\displaystyle \sum_{k=1}^{n}a_{ik}b_{kj}$.

Finally, note that if $A=\{a_{ik}\} \in Mat_{n}(R)$, then $\displaystyle \sum_{i=1}^{n} x_{i}a_{ik} \in F$ and $\hat{f}:x_{j} \mapsto \displaystyle \sum_{i=1}^{n} x_{i}a_{ik}$ is an $R$-homomorphism from $\{x_{i}\}_{i=1}^{n}$ into $F$.
This can be extended to an endomorphism $f \in F$. 
It readily follows that $\psi:Mat_{n}(R) \rightarrow S$ defined by $\{a_{ik}\} \mapsto f$ is a well-defined ring homomorphism. 
Moreover, $\varphi\psi(\{a_{ik}\})=\varphi(f)=\{a_{ik}\}$ and $\psi\varphi(f)=\psi(\{a_{ik}\})=f$.
Thus, $\varphi$ and $\psi$ are inverses, and therefore $\varphi$ is an isomorphism between $S=End_{R}(F)$ and $Mat_{n}(R)$. 
\end{proof}
\end{lemma}

\begin{theorem}\cite{Chatters}\label{semihered iff Mat pp} A ring $R$ is right semi-hereditary if and only if $Mat_{n}(R)$ is a right p.p.-ring for every $0 < n < \omega$.
\begin{proof}
Suppose $R$ is right semi-hereditary.
For $0<n<\omega$, let $F$ be a finitely generated free right $R$-module with basis $\{x_{i}\}_{i=1}^{n}$.
By \Cref{End iso Mat}, $Mat_{n}(R) \cong End_{R}(F)$. 
Therefore, it suffices to show that $S=End_{R}(F)$ is a right p.p.-ring.
Take $s \in S$.
Since $F$ is finitely-generated, $sF$ is a finitely generated submodule of $F$.
Free modules are projective, and thus sF is projective by \Cref{Rot 4.30}.
Since $sF$ is an epimorphic image of $F$, \Cref{projective} shows that $F \cong \ker{s} \bigoplus N$ for some right $R$-module $N$.
Thus, $\ker{s} = eF$ for some nonzero idempotent $e \in S$.
Suppose $r \in ann_{r}(s)= \{t \in S$ $|$ $st(f)=0$ for every $f \in F \}$.
Then, $sr = 0$ and $r \in \ker{s} = eF \subseteq eS$.
On the other hand, suppose $et \in eS$.
Since $sef = 0$ for every $f \in F$, $set(f)=0$ for every $f \in F$.
Hence, $et \in ann_{r}(s)$.
Therefore, $ann_{r}(s)=eS$ and $S=End_{R}(F) \cong Mat_{n}(R)$ is a right p.p.-ring.  

Suppose $Mat_{n}(R)$ is a right p.p.-ring for every $0<n<\omega$.
Let $I$ be a finitely generated right ideal of $R$ with generating set $\{a_{1}, a_{2},...,a_{k}\}$, and take $F$ to be a free right $R$-module with basis $\{x_{1},x_{2},...,x_{k}\}$.
Note that there exists a submodule $K$ of $F$ which is isomorphic to $I$.
Hence, $K$ is also generated by $k$ elements, say $b_{1}, b_{2},...,b_{k}$.
Let $S=Mat_{k}(R) \cong End_{R}(F)$.
For any $f \in F$, there exists $r_{1},r_{2},...,r_{k} \in R$ such that $f = x_{1}r_{1}+x_{2}r_{2}+...+x_{k}r_{k}$.
Let $s \in S$ be the well-defined homomorphism defined by $s(f)=s(x_{1}r_{1}+x_{2}r_{2}+...+x_{n}r_{n}) = b_{1}r_{1}+b_{2}r_{2}+...+b_{n}r_{k}$.
Note that $im (s) = K$ and thus $s:F \rightarrow K$ is an epimorphism.

It will now be shown that $\ker{(s)}=ann_{r}(s)F$.
Here, as before, $ann_{r}(s)$ refers to the annihilator in $S$.
If $y=\displaystyle \sum_{i=1}^{n} t_{i}f_{i} \in ann_{r}(s)F$, then $st_{i}f_{i}=0$ for every $i=1,2,...,n$.
Hence, $y \in \ker{(s)}$.
On the other hand, let $f \in \ker{(s)}$.
Now, $fR$ is a submodule of $F$, and so we can find some $t \in S$ such that $t:F \rightarrow fR$ is an epimorphism and $tf=f$.
Then, for any $x \in F$, $s[t(x)]=s(fr)$ for some $r \in R$.
However, $s(fr)=(sf)r=0$.
Thus, $t \in ann_{r}(s)$ and $f=tf \in ann_{r}(s)F$.
Therefore, $\ker{(s)}=ann_{r}(s)F$.
Moreover, since $Mat_{k}(R) \cong End_{R}(F)$ is a right p.p.-ring by assumption, $ann_{r}(s)=eS$ for some idempotent $e \in S$.
Observe that $SF=F$ since $\displaystyle \sum_{i=1}^{n} s_{i}f_{i} \in F$ for $s_{i} \in S$ and $f_{i} \in F$, and $f = 1_{F}(f) \in SF$ for any $f \in F$.
Hence, $\ker{(s)}=ann_{r}(s)F=eSF=eF$.
Thus, $\ker{(s)}$ is a direct summand of $F$. 
It then follows from \Cref{projective} that $I \cong K$ is projective since $s:F \rightarrow K$ is a an epimorphism.
Therefore, $R$ is a right semi-hereditary ring.
\end{proof}
\end{theorem}

Two idempotents $e$ and $f$ are called {\it orthogonal} if $ef=0$ and $fe=0$.  
If $R$ contains only finite sets of orthogonal idempotents, then being a p.p.-ring is right-left-symmetric.
Moreover, if $R$ is a right (or left) p.p.-ring not containing an infinite set of orthogonal idempotents, then it satisfies both the ascending and descending chain conditions on annihilators (\Cref{RwCC 8.4}).
A ring $R$ satisfies the {\it ascending chain condition} on annihilators if given any ascending chain $I_{0} \subseteq I_{1} \subseteq ... \subseteq I_{n} \subseteq ... $ of annihilators, there exists some $k < \omega$ such that $I_{n}=I_{k}$ for every $n \geq k$.
Similarly, $R$ satisfies the {\it descending chain condition} on annihilators if every descending chain of annihilators terminates for some $k < \omega$. 
Before proving \Cref{RwCC 8.4}, we look at some basic results regarding annihilators and the chain conditions.

\begin{lemma}\label{ann subset ann}Let $S$ and $T$ be subsets of a ring R such that $S \subseteq T$.  Then, $ann_{r}(T) \subseteq ann_{r}(S)$ and $ann_{l}(T) \subseteq ann_{l}(S)$.
\begin{proof}
For $r \in ann_{r}(T)$ and $t \in T$, $tr = 0$.
Let $s \in S \subseteq T$.  Then, $sr = 0$ and hence $r \in ann_{r}(S)$.
Thus, $ann_{r}(T) \subseteq ann_{r}(S)$.
A similar computation shows the theorem holds for left annihilators.
\end{proof}
\end{lemma} 

\begin{lemma}\label{ann of ann}Let U be a subset of a ring R, and let $A = ann_{r}(U) = \{r \in R$ $|$ $ur = 0$ for every $u \in U\}$.
Then, $ann_{r}(ann_{l}(A)) = A$.
\begin{proof}Suppose $r \in ann_{r}(ann_{l}(A))$, and let $u \in U$.
Then, $ua = 0$ for every $a \in A$.
Hence, $u \in ann_{l}(A)$, and thus $ur = 0$.
Therefore, $ann_{r}(ann_{l}(A)) \subseteq A$.
Conversely, suppose $a \in A$.  Then, $ba = 0$ for every $b \in ann_{l}(A)$.
Hence, $a \in ann_{r}(ann_{l}(A))$.  
Therefore, $A \subseteq ann_{r}(ann_{l}(A))$.
\end{proof}
\end{lemma}

\begin{lemma}\label{ann left rt}R satisfies the ascending chain condition on right annihilators if and only if R satisfies the descending chain condition on left annihilators.
\begin{proof}Suppose $R$ satisfies the ascending chain condition on right annihilators.
Let $ann_{l}(U_{1})$ \\$\supseteq ann_{l}(U_{2}) \supseteq ...$ be a descending chain of left annihilators.  
Note that if $ann_{l}(U_{i}) \supseteq ann_{l}(U_{j})$, then  $ann_{r}(ann_{l}(U_{1})) \subseteq ann_{r}(ann_{l}(U_{2})) \subseteq ...$ is an ascending chain of right annihilators by \Cref{ann subset ann}.
By the ascending chain condition on right annihilators, there is some $k < \omega$ such that $ann_{r}(ann_{l}(U_{n})) = ann_{r}(ann_{l}(U_{k}))$ for every $n \geq k$.
Therefore, $ann_{l}(ann_{r}(ann_{l}(U_{n}))) = ann_{l}(ann_{r}(ann_{l}(U_{k})))$ for every $n \geq k$, and by a symmetric version of \Cref{ann of ann} it follows that $ann_{l}(U_{n}) = ann_{l}(U_{n})$ for every $n \geq k$. 
A similar argument shows that the descending chain condition on left annihilators implies the ascending chain condition for right annihilators.
\end{proof}
\end{lemma}

\begin{theorem}\cite{Chatters}\label{RwCC 8.4} Let R be a right p.p.-ring which does not contain an infinite set of orthogonal idempotents.
Then R is also a left p.p.-ring, every right or left annihilator in R is generated by an idempotent, and R satisfies both the ascending and descending chain condition for right annihilators.

\begin{proof}
Let $A = ann_{r}(U)$ for some subset U of R and consider B = $ann_{l}(A)$.
Suppose B contains nonzero orthogonal idempotents $e_{1}, ... , e_{n}$, and let e = $e_{1} + ... + e_{n}$.
Note that e is also an idempotent since $e^{2} = (e_{1} + ... + e_{n})(e_{1} + ... + e_{n}) = e_{1}^{2} + ... + e_{n}^{2} + e_{1}e_{2} + ... + e_{n-1}e_{n} = e_{1} + ... + e_{n} = e$.
Suppose $B = Re$.  The claim is that $A = (1-e)R$, and hence $A$ is generated by an idempotent.
To see this, first note that $ann_{r}(B) = ann_{r}(ann_{l}(A)) = A$ by \Cref{ann of ann}.
Thus, it needs to be shown that $ann_{r}(B) = (1-e)R$.
If $b \in B=Re$, then $b=se$ for some $s\in R$.
For all $r \in R$, we obtain $b(1-e)r = se(1-e)r = (se - se^{2})r = (se - se)r = 0$.
Hence, $(1-e)R \subseteq ann_{r}(B)$.
On the other hand, suppose $r \in ann_{r}(B)$.
Then, $r = r - er + er = (1-e)r + er$.
Note that $e \in B = ann_{l}(A)$, and so $er = 0$ since $r \in ann_{r}(B) = A$.
Thus, $r = (1-e)r \in (1-e)R$, and hence $ann_{r}(B) \subseteq (1-e)R$.
Therefore, if $B=Re$, then $A$ is generated by an idempotent.

If $B \neq Re$, then select $b \in B\char`\\ Re$, and observe $ba = 0$ for every $a \in A$ since $b \neq re$ for any $r \in R$.
Therefore, $B \neq Be$, which implies $B(1-e) \neq 0$.  
Let $0 \neq y \in B(1-e)$, say $y = s(1-e)$ for some $s \in B$.
Since R is a right p.p.-ring, $ann_{r}(y) = (1-f)R$ for some idempotent $f \in R$.  %Why 1 - f instead of f? %
Observe that $f$ is nonzero.
For otherwise, $ann_{r}(y)=R$ and $y=0$, which is a contradiction.
If $0 \neq a \in A$, then $ya = s(1-e)a = sa - sea = 0 - s \cdot 0 = 0$.
Thus, $a \in ann_{r}(y) = (1-f)R$, and so $A \subseteq (1-f)R$.
Hence, $fA \subseteq f(1-f)R=0$ and $f \in ann_{l}(A) = B$.
Observe that $e \in ann_{r}(y) = (1-f)R$ since $ye=s(1-e)e=0$, and so $e = (1-f)t$ for some $t \in R$.
Thus, $(1-f)e = (1-f)(1-f)t = (1-f)t = e$, and so $fe = f(1-f)t = (f - f^{2})t = 0$.
Note also that $fe_{i} = 0$ for $i = 1, ... , n$, since $ye_{i} = s(1-e)e_{i} = s(e_{i} - ee_{i}) = s(e_{i} - e_{i}) = 0$ and hence $e_{i} \in ann_{r}(y)$.  

Let $e_{n+1} = (1-e)f = f - ef$.
Note $e_{n+1}$ is an idempotent since $fe = 0$ and thus $(f - ef)(f - ef) = f - fef - ef + efef = f - 0 - ef + 0 = f - ef$.
Consider $e_{i}$ for some $i = 1, ... , n$.
Then, $e_{n+1}e_{i} = (1-e)fe_{i} = (1-e) \cdot 0 = 0$, and $e_{i}e_{n+1} = e_{i}(1-e)f = (e_{i} - e_{i}e)f = (e_{i} - e_{i})f = 0 \cdot f = 0$.
Thus, $e_{n+1}$ is orthogonal to $e_{1}, ... , e_{n}$.
Furthermore, $e_{n+1}$ is nonzero, since otherwise we have $f = ef$.
This would imply $f = f^{2} = efef = e \cdot 0 \cdot f = 0$, which is a contradiction.  Note also that $e_{n+1} \in B$ since both $e$ and $f$ are in $B$.

Then, $e_{1}, ..., e_{n}, e_{n+1}$ are nonzero orthogonal idempotents contained in $B$.  
As before, if $e = e_{1} + ... + e_{n+1}$ and $B \neq Re$, then there is a nonzero idempotent $e_{n+2} \in B$ orthogonal to $e_{1}, ...,e_{n+1}$.
Since R does not contain any infinite set of orthogonal idempotents, this process must stop for $e_{1}, ..., e_{k}$.  
Thus, for $e = e_{1} + ... + e_{k}$, $B = Re$ and $A = (1-e)R$. 
Therefore, each right and left annihilator is generated by an idempotent.
From a symmetric version of \Cref{pp idempt}, it follows that R is a left p.p.-ring.

Finally, it needs to be shown that R satisfies the ascending and descending chain conditions for right annihilators.
Let $C \subseteq D$ be right annihilators.
Then, there are idempotents $e$ and $f$ such that $C = eR$ and $D = fR$.
Hence, $eR \subseteq fR$, and it follows that $e = fe$. %why? %
Thus, $g = f - ef$ is a nonzero idempotent.
Furthermore, $g$ and $e$ are orthogonal, since $eg = e(f-ef) = ef - e^{2}f = ef - ef = 0$ and $ge = (f - ef)e = fe - efe = e - e^{2} = 0$.
Note that $fR = eR + gR$.
For, if $er + gs \in eR + gR$, then $er + gs = er + (f - ef)s = er + fs + efs \in fR$, and conversely, if $fr \in fR$, then $fr = (f + ef - ef)r = efr + (f - ef)r = efr - gr \in eR + gR$.

Let $I_{1} \subseteq I_{2} \subseteq ... $ be a chain of right annihilators.
Then, for $I_{1} \subseteq I_{2}$, there are idempotents $e$ and $f$ such that $I_{1} = eR$ and $I_{2} = fR$, and there is an idempotent $g$ orthogonal to $e$ such that $I_{2} = I_{1} + gR$.
It then follows that $I_{3} = I_{1} + gR + hR$ for some idempotent $h$ orthogonal to both $e$ and $g$.  
Since $R$ does not contain an infinite set of orthogonal idempotents, this must terminate with some $k < \omega$ so that $I_{n} = I_{k}$ for every $n \geq k$.  Therefore, $R$ satisfies the ascending chain condition on right annihilators.  The descending chain condition on right annihilators follows from \Cref{ann left rt}.  
\end{proof}
\end{theorem}

\chapter{Homological Algebra}

Before discussing torsion-freeness and non-singularity of modules, we need some basic results in Homological Algebra regarding tensor products, flat modules, and functors.

\section{Tensor Products}

Let $A$ be a right $R$-module, $B$ a left $R$-module, and $G$ any Abelian group.
A function $f: A \times B \rightarrow G$ is called {\it$R$-biadditive}, or {\it $R$-bilinear}, if the following conditions are satisfied:
\begin{enumerate}[(i)]
\item For each $a, a^{\prime} \in A$ and $b \in B$, $f(a+a^{\prime},b)=f(a,b)+f(a^{\prime},b)$,
\item For each $a \in A$ and $b, b^{\prime} \in B$, $f(a,b+b^{\prime})=f(a,b)+f(a,b^{\prime})$, 
\item For each $a \in A$, $b \in B$, and $r \in R$, $f(ar,b)=f(a,rb)$.
\end{enumerate}
Note that in general $f(a+a^{\prime},b+b^{\prime}) \neq f(a,b)+f(a^{\prime},b^{\prime})$.
The {\it tensor product} of $A$ and $B$, denoted $A \bigotimes_{R} B$, is an Abelian group and an $R$-biadditive function $h:A \times B \rightarrow A \bigotimes_{R} B$ having the universal property that whenever $G$ is an Abelian group and $g:A \times B \rightarrow G$ is $R$-biadditive, there is a unique map $f:A \bigotimes_{R} B \rightarrow G$ such that $g=fh$.

\begin{prop}\cite{Rotman}\label{tensor exist}
Let $R$ be a ring.  Given a right $R$-module $A$ and a left $R$-module $B$, the tensor product $A \bigotimes_{R} B$ exists.
\begin{proof}
Let $F$ be a free Abelian group with basis $A \times B$, and let $U$ be a subgroup of $F$ generated by all elements of the form $(a+a^{\prime},b)-(a,b)-(a^{\prime},b)$, $(a,b+b^{\prime})-(a,b)-(a,b^{\prime})$, or $(ar,b)-(a,rb)$, where $a, a^{\prime} \in A$, $b, b^{\prime} \in B$, and $r \in R$.
Define $A \bigotimes_{R} B$ to be $F \slash U$, and denote $(a,b)+U \in F \slash U$ as $a \otimes b$.
In addition, let $h:A \times B \rightarrow A \bigotimes_{R} B$ be defined by $(a,b) \mapsto a \otimes b$.
Observe that $h$ is a well-defined $R$-biadditive map. 
For if $a, a^{\prime} \in A$ and $b \in B$, then $h(a+a^{\prime},b)=(a+a^{\prime},b)+U=(a+a^{\prime},b)-[(a+a^{\prime},b)-(a,b)-(a^{\prime},b)]+U=[(a,b)+U] + [(a^{\prime},b)+U]=h(a,b)+h(a^{\prime},b)$.
Similarly, $h(a,b+b^{\prime})=h(a,b)+h(a,b^{\prime})$ for $b,b^{\prime} \in B$, and $h(ar,b)=(ar,b)+U=(ar,b)-[(ar,b)-(a,rb)]+U=(a,rb)+U=h(a,rb)$ for $r \in R$.

Let $G$ be any Abelian group and $g:A \times B \rightarrow G$ any $R$-biadditive map.
For $F \slash U$ to be a tensor product, it needs to be shown that there is a function $\varphi:A \bigotimes_{R} B =F \slash U \rightarrow G$ such that $g =\varphi h$.
Define $\hat{f}:A \times B \rightarrow G$ by $(a,b) \mapsto g(a,b)$.
Each element of $F$ is of the form $\sum_{A \times B}$ $(a,b)n_{(a,b)}$, where $n_{(a,b)}=0$ for all but finitely many $(a,b) \in A \times B$.
Let $f$ be defined by $\sum_{A \times B}$ $(a,b)n_{(a,b)} \mapsto \sum_{A \times B}$ $\hat{f}[(a,b)]n_{(a,b)}$.
This is clearly well-defined since $\hat{f}$ is well-defined.
Moreover, $f[(a,b)]= \hat{f}[(a,b)]$ for $(a,b) \in A \times B$, and thus $f$ extends $\hat{f}$ to a function on $F$.
Note that if $k$ is another extension of $\hat{f}$, then $k$ must equal $f$ since they are equal on the generating set $A \times B$.
Hence, $f$ is a unique extension.
Also observe that $f$ is a homomorphism since, given $x,y \in F$, $f(x+y)=f(\sum_{A \times B}$ $(a,b)n_{(a,b)}+\sum_{A \times B}$ $(a^{\prime},b^{\prime})m_{(a,b)})
\\=\sum_{A \times B}$ $\hat{f}[(a,b)]n_{(a,b)}+\sum_{A \times B}$ $\hat{f}[(a^{\prime},b^{\prime})]m_{(a,b)}=f(x)+f(y)$.

It readily follows from $g$ being $R$-biadditive that the homomorphism $f:F \rightarrow G$ which we have just constructed is also $R$-biadditive.
To see this, observe that if $a, a^{\prime} \in A$ and $b \in B$, then $f[(a+a^{\prime},b)]-f[(a,b)]-f[(a^{\prime},b)]=g[(a+a^{\prime},b)]-g[(a,b)]-g[(a^{\prime},b)]=0$.
The other two conditions are satisfied with similar computation.
Thus, we have that $f(U)=0$.
Define $\varphi:F \slash U = A \bigotimes_{R} B \rightarrow G$ by $\varphi(x+U)=f(x)$.
If $x+U=x^{\prime}+U$, then $x-x^{\prime} \in U$ and hence $f(x-x^{\prime}) \in f(U)=0$.
Thus, $f(x)=f(x^{\prime})$ and $\varphi$ is well-defined.
Furthermore, $\varphi h(a,b)=\varphi [a \otimes b]=\varphi[(a,b)+U]=f[(a,b)]=g[(a,b)]$.
Therefore $A \bigotimes_{R} B = F \slash U$ is a tensor product.
\end{proof}
\end{prop}

\begin{prop}
Let $R$ be a ring, $A$ a right $R$-module, and $B$ a left $R$-module.
Then, the tensor product $A \bigotimes_{R} B$ is unique up to isomorphism.
\begin{proof}
It has already been shown that $A \bigotimes_{R} B$ exists.
Suppose $H$ and $H^{\prime}$ are both tensor products, and let $h:A \times B \rightarrow H$ and $h^{\prime}:A \times B \rightarrow H^{\prime}$ be the respective $R$-biadditive functions having the universal property. 
Then, there exists a function $f:H \rightarrow H^{\prime}$ such that $h^{\prime}=fh$ and a function $f^{\prime}:H^{\prime} \rightarrow H$ such that $h=f^{\prime}h^{\prime}$.
Hence, $h=f^{\prime}fh$ and $h^{\prime}=ff^{\prime}h^{\prime}$.
That is, $f^{\prime}f \cong 1_{H}$ and $ff^{\prime} \cong 1_{H^{\prime}}$.
Therefore, $f:H \rightarrow H^{\prime}$ is an isomorphism.
\end{proof}
\end{prop}

Each element of $A \bigotimes_{R} B$ is a finite sum of the form $\displaystyle{\sum_{i=1}^{n} (a_{i} \otimes b_{i})}$.
The elements $a \otimes b$ that generate $A \bigotimes_{R} B$ are referred to as {\it tensors}.
Given $a, a^{\prime} \in A$, $b, b^{\prime} \in B$, and $r \in R$, the following properties hold for tensors:
\begin{enumerate}[(i)]
\item $(a+a^{\prime}) \otimes b=a \otimes b + a^{\prime} \otimes b$,
\item $a \otimes (b+b^{\prime})=a \otimes b + a \otimes b^{\prime}$,
\item $ar \otimes b = a \otimes rb$.
\end{enumerate}
These properties can be proved in a method similar to that used in the proof of \Cref{tensor exist} to show that $h:A \times B \rightarrow A \bigotimes_{R} B$ defined by $(a,b) \mapsto a \otimes b$ is $R$-biadditive.

\begin{prop}\cite{Rotman}\label{induced tensor}
Let $R$ be a ring, $A, A^{\prime} \in Mod_{R}$, and $B, B^{\prime} \in$  $_{R}Mod$.
If $f:A \rightarrow A^{\prime}$ and $g:B \rightarrow B^{\prime}$ are $R$-homomorphisms, then there is an induced map $f \otimes g: A \bigotimes_{R} B \rightarrow A^{\prime} \bigotimes_{R} B^{\prime}$ such that $(f \otimes g)(a \otimes b)=f(a) \otimes g(b)$.
\begin{proof}
Let $h:A \times B \rightarrow A \bigotimes_{R} B$ and $h^{\prime}:A^{\prime} \times B^{\prime} \rightarrow A^{\prime} \bigotimes_{R} B^{\prime}$ be the respective $R$-biadditive maps with the universal tensor property.
Define $\varphi:A \times B \rightarrow A^{\prime} \times B^{\prime}$ by $\varphi(a,b)=(f(a),g(b))$.
It then follows that $h^{\prime}\varphi:A \times B \rightarrow A^{\prime} \bigotimes_{R} B^{\prime}$ is $R$-biadditive.
For if $a, a^{\prime} \in A$ and $b \in B$, then $h^{\prime}\varphi(a+a^{\prime},b)=h^{\prime}(f(a+a^{\prime}),g(b))=h^{\prime}[f(a)+f(a^{\prime}),g(b)]=h^{\prime}[f(a),g(b)]+h^{\prime}[f(a^{\prime}),g(b)]=h^{\prime}\varphi(a,b)+h^{\prime}\varphi(a^{\prime},b)$.
Similarly, $h^{\prime}\varphi(a,b+b^{\prime})=h^{\prime}\varphi(a,b)+h^{\prime}\varphi(a,b^{\prime})$ and $h^{\prime}\varphi(ar,b)=h^{\prime}\varphi(a,rb)$ for $b^{\prime} \in B$ and $r \in R$.
By the universal property of the $R$-biadditive map $h$, there exists a map $\hat{\varphi}:A \bigotimes_{R} B \rightarrow A^{\prime} \bigotimes_{R} B^{\prime}$ such that $h^{\prime}\varphi = \hat{\varphi}h$.
Hence, $\hat{\varphi}(a \otimes b)=\hat{\varphi}h(a,b)=h^{\prime}\varphi(a,b)=h^{\prime}[f(a),g(b)]=f(a) \otimes g(b)$.
Therefore, $f \otimes g=\hat{\varphi}$ is an induced map satisfying $(f \otimes g)(a \otimes b)=f(a) \otimes g(b)$.
\end{proof}
\end{prop}

The following lemmas will be needed in a later section:

\begin{lemma}\cite{Fuchs}\label{tensor zero}
Let $R$ be a ring, $A$ a right $R$-module, and $B$ a left $R$-module.
If $a \otimes b$ is a tensor in $A \bigotimes_{R} B$, then $a \otimes b=0$ if and only if there exists $a_{1},a_{2},...,a_{k} \in A$ and $r_{1},r_{2},...,r_{k} \in R$ such that $a=a_{1}r_{1}+a_{2}r_{2}+...+a_{k}r_{k}$ and $r_{j}b=0$ for $j=1,2,...,k$.
\end{lemma}

\begin{lemma}\label{Tensor Nat Iso}For a left R-module M, there is an R-module isomorphism \\$\varphi$ : R $\bigotimes_{R} M \rightarrow M$ given by $\varphi (r \otimes m) = rm$.  Here, $R$ is viewed as a right R-module.
Similarly, $N \bigotimes_{R} R \cong N$ for a right $R$-module $N$.
\begin{proof}First, observe that $R \times M \xrightarrow{\psi} M$ given by $\psi ((r, m)) = rm$ is R-biadditive.
%For $\psi (r + r^{'}, m) = (r + r^{'})m = rm + r^{'}m = \psi(r, m) + \psi(r^{'}, m)$, $\psi (r, m + m^{'}) = r(m + m^{'}) = rm + r m^{'} = \psi(r, m) + \psi(r, m^{'})$, and given any $s \in R$, $\psi (rs, m) = (rs)m = r(sm) = \psi(r, sm)$.
Thus, we can define an R-module homomorphism $R \bigotimes_{R} M \xrightarrow{\varphi} M$ that sends each $r \otimes m \in R \bigotimes_{R} M$ to $rm$.
In other words, $\varphi(r \otimes m) = \psi(r, m)$.
Note that for every $s \in R$, $\varphi(s(r \otimes m)) = \varphi(sr \otimes m) = (sr)m = s(rm) = s\varphi(r \otimes m)$.

Let $\alpha : M \rightarrow R \bigotimes_{R} M$  be defined by $\alpha(m) = 1 \otimes m$.
Clearly $\alpha$ is a well-defined R-module homomorphism since $\alpha(m + n) = 1 \otimes (m +n) = 1 \otimes m + 1 \otimes n = \alpha(m) + \alpha(n)$, and $\alpha(rm) = 1 \otimes rm = 1r \otimes m = 1 \otimes m$.
It follows that $\alpha\varphi(r \otimes m) = \alpha(rm) = 1 \otimes rm = 1r \otimes m = r \otimes m$, and
$\varphi\alpha(m) = \varphi(1 \otimes m) = 1m = m$.
Thus, $\varphi$ is a bijection and hence an R-module isomorphism.
\end{proof}
\end{lemma}

\begin{lemma}\cite{Rotman}\label{Right Exactness}If $A \xrightarrow{i} B \xrightarrow{p} C \rightarrow 0$ is an exact sequence of left R-modules, 
then for any right R-module M, $M \bigotimes_{R} A \xrightarrow{1 \otimes i} M \bigotimes_{R} B \xrightarrow{1 \otimes p} M \bigotimes_{R} C \rightarrow 0$ is an exact sequence.
\begin{proof}
For $M \bigotimes_{R} A \xrightarrow{1 \otimes i} M \bigotimes_{R} B \xrightarrow{1 \otimes p} M \bigotimes_{R} C \rightarrow 0$ to be exact, it needs to be shown that $im (1 \otimes i) = \ker{(1 \otimes p)}$ and $1 \otimes p$ is surjective.
Since $im (i) = \ker{(p)}$ and hence $pia = 0$ for every $a \in A$, it readily follows that $im (1 \otimes i) \subseteq \ker{(1 \otimes p)}$.
For if $\sum (m_{j} \otimes a_{j}) \in M \bigotimes_{R} A$, then $(1 \otimes p)(1 \otimes i)[\sum (m_{j} \otimes a_{j})] = (1 \otimes p)[\sum (1 \otimes i)(m_{j} \otimes a_{j})] = (1 \otimes p)[\sum (m_{j} \otimes ia_{j})] = \sum (1 \otimes p)(m_{j} \otimes ia_{j})=\sum (m_{j} \otimes pia_{j}) = \sum (m_{j} \otimes 0) = 0$.
To see that $im(1 \otimes i) = \ker{(1 \otimes p)}$, first note that since $im (1 \otimes i)$ is contained in the kernel of $1 \otimes p$, there is a uniqe homomorphism $\varphi: M \bigotimes_{R} B \slash im(1 \otimes i) \rightarrow M \bigotimes_{R} C$ such that $\varphi[(m \otimes b) + im(1 \otimes i)] = (1 \otimes p)(m \otimes b) = m \otimes pb$ \cite[Ch. IV, Theorem 1.7]{Hungerford}.

It can be shown that $\varphi$ is an isomorphism, and from this it will follow that $im(1 \otimes i) = \ker{(1 \otimes p)}$.
Note that since the sequence $A \xrightarrow{i} B \xrightarrow{p} C \rightarrow 0$ is exact and hence $p$ is surjective, for every $c \in C$ there exists an element $b \in B$ such that $pb = c$.
Let the function $f:M \times C \rightarrow M \bigotimes_{R} B \slash im(1 \otimes i)$ be defined by $(m,c) \mapsto p \otimes b$.
If there is another element $b_{0} \in B$ such that $pb_{0} = c$, then $p(b - b_{0})=pb-pb{0}=c-c=0$.
Hence, $b-b_{0} \in \ker{(p)}=im(i)$.
Thus, there is an $a \in A$ such that $ia = b - b_{0}$, and it then follows that $m \otimes b - m \otimes b_{0} = m \otimes (b-b_{0}) = m \otimes ia \in im(1 \otimes i)$.
Hence, $(m \otimes b - m \otimes b_{0}) + im(1 \otimes i) = 0$, and therefore $f$ is well-defined.
Furthermore, it is easily seen that $f$ is an $R$-biadditive function.
Thus, if $h: (m,c) \mapsto m \otimes c$ is the biadditive function of the tensor product, then there is a homomorphism $\psi: M \bigotimes_{R} C \rightarrow M \bigotimes_{R} B \slash im(1 \otimes i)$ such that $\psi h = f$.
In other words, $\psi(m \otimes c) = (m \otimes b)+im(1 \otimes i)$.

Observe that $\psi\varphi[(m \otimes b)+im(1 \otimes i)] = \psi(m \otimes pb) = \psi(m \otimes c) = (m \otimes b) + im(1 \otimes i)$ and $\varphi\psi(m \otimes c)=\varphi[(m \otimes b)+im(1 \otimes i)]=m \otimes pb= m \otimes c$.
Thus, $\varphi$ is an isomorphism with inverse $\psi$.
Now, let $\pi:M \bigotimes_{R} B \rightarrow M \bigotimes_{R} B \slash im(1 \otimes i)$ be the canonical epimorphism given by $m \otimes b \mapsto m \otimes b + im(1 \otimes i)$.
Then, $\varphi\pi(m \otimes b) = \varphi[(m \otimes b) + im(1 \otimes i)] = m \otimes pb = (1 \otimes p)(m \otimes b).$
Hence, $\varphi\pi = 1 \otimes p$.
Therefore, since $\varphi$ is an isomorphism, $\ker{(1 \otimes p)}=\ker{(\varphi\pi)}= \ker{(\pi)}=im(1+i)$.

Finally, it needs to be shown that $1 \otimes p$ is surjective.
Let $\sum (m_{j} \otimes c_{j}) \in M \bigotimes_{R} C$. 
Since $p$ is surjective, for each $j$, there exists an element $b_{j} \in B$ such that $pb_{j} = c_{j}$.
Thus, $(1 \otimes p)[\sum (m_{j} \otimes b_{j})] = \sum (1 \otimes p)(m_{j} \otimes b_{j}) = \sum (m_{j} \otimes pb_{j}) = \sum (m_{j} \otimes c_{j})$.
Therefore, $1 \otimes p$ is surjective and the sequence $M \bigotimes_{R} A \xrightarrow{1 \otimes i} M \bigotimes_{R} B \xrightarrow{1 \otimes p} M \bigotimes_{R} C \rightarrow 0$ is exact.
\end{proof}
\end{lemma}

A right R-module M is {\it flat} if $0 \rightarrow M \bigotimes_{R} A \xrightarrow{1_{M} \otimes \varphi} M \bigotimes_{R} B \xrightarrow{1_{M} \otimes \psi} M \bigotimes_{R} C \rightarrow 0$ is an exact sequence of Abelian groups whenever
$0 \rightarrow A \xrightarrow{\varphi} B \xrightarrow{\psi} C \rightarrow 0$ is an exact sequence of left R-modules.

\begin{prop}\cite{Rotman}\label{direct sum of flat}
Let R be a ring and let $\{M_{i}\}_{i \in I}$ be a collection of right R-modules for some index set I.  Then, the direct sum $\bigoplus_{I} M_{i}$ is flat if and only if $M_{i}$ is flat for every $i \in I$.  Moreover, $R$ is flat as a right R-module, and any projective right R-module P is flat.
\begin{proof}
First note that if $0 \rightarrow A \xrightarrow{\varphi} B \xrightarrow{\psi} C \rightarrow 0$ is an exact sequence of left $R$-modules, then $M \bigotimes_{R} A \xrightarrow{1_{M} \otimes \varphi} M \bigotimes_{R} B \xrightarrow{1_{M} \otimes \psi} M \bigotimes_{R} C \rightarrow 0$ is exact by \Cref{Right Exactness}.
Thus, $M$ is flat if and only if $1_{M} \otimes \varphi$ is a monomorphism whenever $\varphi$ is a monomorphism.

Suppose $A$ and $B$ are left $R$-modules and let $\varphi:A \rightarrow B$ be a monomorphism.
For $\bigoplus_{I} M_{i}$ to be flat, it needs to be shown that $1 \otimes \varphi:(\bigoplus_{I} M_{i}) \bigotimes_{R} A \rightarrow (\bigoplus_{I} M_{i}) \bigotimes_{R} B$ is a monomorphism.
By \cite[Theorem 2.65]{Rotman}, there exist isomorphisms $f:(\bigoplus_{I} M_{i}) \bigotimes_{R} A \rightarrow (\bigoplus_{I} M_{i} \bigotimes_{R} A)$ and $g:(\bigoplus_{I} M_{i}) \bigotimes_{R} B \rightarrow (\bigoplus_{I} M_{i} \bigotimes_{R} B)$ defined by $f:(x_{i}) \otimes a \mapsto (x_{i} \otimes a)$ and $g:(x_{i}) \otimes b \mapsto (x_{i} \otimes b)$.
Furthermore, since $1_{M_{i}} \otimes \varphi$ is a homomorphism for each $i \in I$, there is a homomorphism $\psi:\bigoplus_{I}(M_{j} \bigotimes_{R} A) \rightarrow \bigoplus_{I}(M_{j} \bigotimes_{R} B)$ such that $(x_{i} \otimes a) \mapsto (x_{i} \otimes \varphi(a))$.
Observe that $\psi$ is a monomorphism if and only if $1_{M_{i}} \otimes \varphi$ is a monomorphism for each $i \in I$.
It then follows that $\psi f=g(1 \otimes \varphi)$ since $\psi f[(x_{i}) \otimes a]=\psi(x_{i} \otimes a)=x_{i} \otimes \varphi(a)=g[(x_{i}) \otimes \varphi(a)]=g(1 \otimes \varphi)[(x_{i} \otimes a)]$.
Therefore, $\bigoplus_{I} M_{i}$ is flat if and only if $1 \otimes \varphi$ is a monomorphism if and only if $\psi$ is a monomorphism if and only if $1_{M_{i}} \otimes \varphi$ is a monomorphism for each $i$ if and only if $M_{i}$ is flat for each $i$.

To see that $R$ is flat as a right $R$-module, note that \Cref{Tensor Nat Iso} gives isomorphisms $f:A \rightarrow R \bigotimes_{R} A$ and $g:B \rightarrow R \bigotimes_{R} B$ defined by $f(a)=1_{R} \otimes a$ and $g(b)=1_{R} \otimes b$.
Observe that $(1_{R} \otimes \varphi)f(a)=(1_{R} \otimes \varphi)(1_{R} \otimes a)=1_{R} \otimes \varphi(a)=g\varphi(a)$.
Hence, $(1_{R} \otimes \varphi)=g\varphi f^{-1}$, which is a monomorphism.
Therefore, $R$ is flat as a right $R$-module.

Let $P$ be a projective right $R$-module.
Then there is a free right $R$-module $F$ and an $R$-module $N$ such that $F=P \bigoplus N$.
As a free module, $F$ is a direct sum of copies of $R$, which is flat.
Hence, $F$ is also flat.
Therefore, $P$ is flat as a direct summand of $F$.
\end{proof}
\end{prop}

\section{Bimodules and the Hom and Tensor Functors}

Let $A$ be a right $R$-module.  
Consider the functor $T_{A}:$ $_{R}Mod \rightarrow Ab$ defined by $T_{A}(B)=A \bigotimes_{R} B$ with induced map $T_{A}(\varphi)=1_{A} \otimes \varphi:A \bigotimes_{R} B \rightarrow A \bigotimes_{R} B^{\prime}$, where $Ab$ is the category of all Abelian groups and $\varphi \in Hom_{R}(B,B^{\prime})$ for left $R$-modules $B$ and $B^{\prime}$.
Observe that $T_{A}(\varphi)(a \otimes b)=a \otimes \varphi(b)$.
$T_{A}$ is sometimes denoted $T_{A}(\underline{\hspace{.3cm}})=A \bigotimes_{R} \underline{\hspace{.3cm}}$.
Similarly, the functor $T_{B}(A)=A \bigotimes_{R} B$ with induced map $\psi \otimes 1_{B}$ can be defined for a left $R$-module $B$ and $\psi \in Hom_{R}(A,A^{\prime})$.
We also consider the functor $Hom_{R}(A, \underline{\hspace{.3cm}}):Mod_{R} \rightarrow Ab$ with induced map $f_{*}:Hom_{R}(A,B) \rightarrow Hom_{R}(A,C)$ defined by $f_{*}(h)=fh$, where $f:B \rightarrow C$ is a homomorphism for right $R$-modules $B$ and $C$.

Let $R$ and $S$ be rings and let $M$ be an Abelian group which has both a left $R$-module structure and a right $S$-module structure.
Then, $M$ is an {\it $(R,S)$-bimodule} if $(rx)s=r(xs)$ for every $r \in R$, $s \in S$, and $x \in M$.
This is sometimes denoted $_{R}M_{S}$.
In particular, if $A$ is a right $R$-module and $E=End_{R}(A)$, then $M$ is an $(E,R)$-bimodule.
Note that for $x \in M$ and $\alpha \in E$, scalar multiplication $\alpha x$ is defined as $\alpha(x)$. 

\begin{prop}\label{Hom module}
Let $R$ and $S$ be rings.
Suppose $M$ is an $(R,S)$-bimodule and $N$ is a right $S$-module.
Then, $Hom_{S}(M_{S},N_{S})$ is a right $R$-module and $Hom_{S}(N_{S},M_{S})$ is a left $R$-module.
\begin{proof}
First, observe that $Hom_{S}(M_{S},N_{S})$ is an Abelian group.
For if $f,g \in Hom_{S}(M_{S},N_{S})$, then $f(xr)=f(x)r$ and $g(xr)=g(x)r$ for every $r \in R$.
Hence, $f+g \in Hom_{S}(M_{S},N_{S})$ since $(f+g)(xr)=f(xr)+g(xr)=f(x)r+g(x)r=(f+g)(x)r$.
Moreover, if $h \in Hom_{S}(M_{S},N_{S})$, then $[f+(g+h)](x)=f(x)+(g+h)(x)=f(x)+g(x)+h(x)=(f+g)(x)+h(x)=[(f+g)+h](x)$.
Hence, $Hom_{S}(M_{S},N_{S})$ is associative.
Furthermore, the map $\alpha:a \mapsto 0$ acts as the zero element.
Finally, note that if $f \in Hom_{S}(M_{S},N_{S})$, then $g:M \rightarrow N$ defined by $g(x)=-f(x)$ is such that $(f+g)(x)=f(x)+g(x)=f(x)-f(x)=0$.
Hence, every element of $Hom_{S}(M_{S},N_{S})$ has an inverse.
Therefore, $Hom_{S}(M_{S},N_{S})$ is an Abelian group.

Now, let $\varphi \in Hom_{S}(M_{S},N_{S})$, $r,r^{\prime} \in R$, and $x \in M$.
Define the right $R$-module structure on $Hom_{S}(M_{S},N_{S})$ by $(\varphi r)(x)=\varphi(rx)$.
Then, $(\varphi + \psi)(r)(x)=(\varphi r + \psi r)(x)=(\varphi r)(x)+(\psi r)(x)=\varphi(rx)+\psi(rx)=(\varphi + \psi)(rx)$ for $\psi \in Hom_{S}(M_{S},N_{S})$.
Moreover, $[\varphi(r+r^{\prime})](x)=\varphi[(r+r^{\prime})x]=\varphi[rx+r^{\prime}x]=\varphi(rx)+\varphi(r^{\prime}x)=(\varphi r)(x)+(\varphi r^{\prime})(x)$ for $r^{\prime} \in R$.
Finally, observe that $[\varphi (rr^{\prime})](x)=\varphi[(rr^{\prime})(x)]=\varphi[r(r^{\prime}x)]=(\varphi r)(r^{\prime}x)$.
Therefore, $Hom_{S}(M_{S},N_{S})$ satisfies the conditions of a right $R$-module.
Similarly, $Hom_{S}(N_{S},M_{S})$ is a left $R$-module with $(r\pi)(x)=r\pi(x)$ for any $\pi \in Hom_{S}(N_{S},M_{S})$.
\end{proof}
\end{prop}

\begin{prop}\cite{Rotman}\label{Rot 2.51}
Let $R$ be a subring of $S$.  Suppose $M$ is an $(R,S)$-bimodule and $A$ is a right $R$-module.  
Then, $A \bigotimes_{R} M$ is a right $S$-module.
In particular, $S$ is an $(R,S)$-bimodule and hence $A \bigotimes_{R} S$ is a right $S$-module.
\begin{proof}
Let $y=\displaystyle{\sum_{i=1}^{n} (a_{i} \otimes x_{i})}$ $\in A \bigotimes_{R} M$ and let $s \in S$.
Define the right $S$-module structure on $A \bigotimes_{R} M$ by $(\displaystyle{\sum_{i=1}^{n} (a_{i} \otimes x_{i})})s = \displaystyle{\sum_{i=1}^{n} (a_{i} \otimes x_{i}s)}$.
To see that this does define a right $S$-module, consider the well-defined map $\mu_{s}:M \rightarrow M$ defined by $\mu_{s}(x)=xs$.
By the bimodule structure of $M$, $r\mu_{s}(x)=r(xs)=(rx)s=\mu_{s}(rx)$ for $r \in R$.
Hence, $\mu_{s} \in Hom_{R}(M,M)$.
Consider the functor $T_{A}(\underline{\hspace{.3cm}})=A \bigotimes_{S} \underline{\hspace{.3cm}}$.
By \Cref{induced tensor}, there is a well-defined homomorphism $T_{A}(\mu_{s})=1_{A} \otimes \mu_{s}:A \bigotimes_{R} M \rightarrow A \bigotimes_{R} M$ such that $(1_{A} \otimes \mu_{s})(a \otimes x)=a \otimes \mu_{s}(x)=a \otimes xs$.
If the element $ys$ is defined by $ys=(1_{A} \otimes \mu_{s})(y)=(1_{A} \otimes \mu_{s})(\displaystyle{\sum_{i=1}^{n} (a_{i} \otimes x_{i})})=\displaystyle{\sum_{i=1}^{n} (1_{A} \otimes \mu_{s})(a_{i} \otimes x_{i})}=\displaystyle{\sum_{i=1}^{n} (a_{i} \otimes x_{i}s)}$, then the $S$-module structure is well-defined since $(1_{A} \otimes \mu_{s})$ is a well-defined homomorphism and $\displaystyle{\sum_{i=1}^{n} (a_{i} \otimes x_{i}s)}$ $\in A \bigotimes_{R} M$.
The remaining right $S$-module conditions follow readily.
Moreover, it is easy to see that $S$ satisfies the conditions of an $(R,S)$-bimodule.
Therefore, given any right $R$-module $A$, $A \bigotimes_{R} S$ is a right $S$-module. 
\end{proof}
\end{prop}

\begin{prop}\label{Functor into ModR ModS}
Let $R \leq S$ be rings and let $M$ be an $(R,S)$-bimodule.
Then, the following hold: 
\begin{enumerate}[(a)]
\item The functor $T_{M}(\underline{\hspace{.3cm}})=\underline{\hspace{.3cm}} \bigotimes_{R} M:Mod_{R} \rightarrow Ab$ is actually a functor $Mod_{R} \rightarrow Mod_{S}$.
\item The functor $Hom_{S}(M,\underline{\hspace{.3cm}}):Mod_{S} \rightarrow Ab$ is actually a functor $Mod_{S} \rightarrow Mod_{R}$.
\end{enumerate}
\begin{proof}
($a$): It has already been shown in \Cref{Rot 2.51} that $T_{M}(A)=A \bigotimes_{R} M$ is a right $S$-module for any right $R$-module $A$.
It needs to be shown that if $\psi \in Hom_{R}(A,A^{\prime})$ for $A^{\prime} \in Mod_{R}$, then $T_{M}(\psi)=\psi \otimes 1_{M} \in Hom_{S}(A \bigotimes_{R} M, A^{\prime} \bigotimes_{R}M)$.
In other words, it needs to be shown that $\psi \otimes 1_{M}$ is an $S$-homomorphism.
Let $s \in S$.
Then, $(\psi \otimes 1_{M})(a \otimes x)s=(\psi(a) \otimes x)s=\psi(a) \otimes xs=(\psi \otimes 1_{M})(a \otimes xs)=(\psi \otimes 1_{M})[(a \otimes x)s]$.
Thus, $T_{M}(\psi)$ is a morphism in $Mod_S$, and therefore $T_{M}(\underline{\hspace{.3cm}})$ is a functor with values in $Mod_{S}$.

($b$): Given any right $S$-module $N$, $Hom_{S}(M,N)$ is a right $R$-module by \Cref{Hom module}.
It needs to be shown that if $f:N \rightarrow N^{\prime}$ is a homomorphism for $N, N^{\prime} \in Mod_{S}$, then the induced map $f_{*}=Hom_{R}(M,f): Hom_{S}(M,N) \rightarrow Hom_{S}(M,N^{\prime})$ defined by $f_{*}(\varphi)=f\varphi$ is an $R$-homomorphism.
Note that if $\varphi, \psi \in Hom_{S}(M,N)$, then $f(\varphi +\psi)=f\varphi+f\psi$.
Hence, $f_{*}$ is a homomorphism since $f_{*}(\varphi +\psi)=f(\varphi+\psi)=f\varphi+f\psi=f_{*}\varphi+f_{*}\psi$.
Let $r \in R$.
Observe that $(\varphi r)(x)=\varphi(rx)$ by \Cref{Hom module}.
Moreover, since $M$ has a left $R$-module structure and $f\varphi$ is an element of the right $R$-module $Hom_{S}(M,N^{\prime})$, \Cref{Hom module} also shows that $[f\varphi(x)]r=f[\varphi r](x)=f\varphi(rx)$ for $x \in M$.
Thus, $[f_{*}(\varphi(x))]r=[f\varphi(x)]r=f\varphi(rx)=f_{*}[\varphi(rx)]=f_{*}[(\varphi r)(x)]$.
Hence, $f_{*}$ is an $R$-homomorphism, and therefore $Hom_{S}(M,\underline{\hspace{.3cm}})$ is a functor with values in $Mod_{R}$.
\end{proof}
\end{prop}

The following lemmas will be used later to show $Mod_{R} \cong Mod_{Mat_{n}}(R)$.  The proofs are omitted and can be found in {\it Rings and Categories of Modules} by Frank Anderson and Kent Fuller.

\begin{lemma}\cite[Proposition 20.10]{Anderson}\label{And 20.10}
Let $R$ and $S$ be rings, $M$ a right $R$-module, $N$ a right $S$-module, and $P$ an $(S,R)$-bimodule.
If $M$ is finitely generated and projective, then $\mu:N \bigotimes_{S} Hom_{R}(M,P) \rightarrow Hom_{R}(M,N \bigotimes_{S} P)$ defined by $\mu(y \otimes f)(x)=y \otimes f(x)$ is a natural isomorphism.
Here, $x \in M$, $y \in N$, and $f \in Hom_{R}(M,P)$.
\end{lemma}

\begin{lemma}\cite[Proposition 20.11]{Anderson}\label{And 20.11}
Let $R$ and $S$ be rings, $M$ a right $R$-module, $N$ a left $S$-module, and $P$ an $(S,R)$-bimodule.
If $M$ is finitely generated and projective, then $\nu:Hom_{R}(P,M) \bigotimes_{S} N \rightarrow Hom_{R}(Hom_{S}(N,P),M)$ defined by $\nu(f \otimes y)(g)=fg(y)$ is a natural isomorphism.
Here, $f \in Hom_{R}(P,M)$, $g \in Hom_{S}(N,P)$, and $y \in N$.
\end{lemma}
\section{The Tor Functor}

Consider the exact sequence $P = \cdot\cdot\cdot \rightarrow P_{2} \xrightarrow{d_{2}} P_{1} \xrightarrow{d_{1}} P_{0} \xrightarrow{\epsilon} A \rightarrow 0$ of right $R$-modules, where $P_{j}$ is projective for every $j$. 
Such an exact sequence is called a {\it projective resolution} of the right $R$-module $A$. 
Note that a projective resolution can be formed for any projective right $R$-module $A$ since every right $R$-module is the epimorphic image of a projective right $R$-module.  
Define the {\it deleted projective resolution}, denoted $P_{A}$, by removing the morphism $\epsilon$ and the right R-module A.  Note that the projective resolution is an exact sequence, and hence $im({d_{i+1}}) = \ker{(d_{i})}$.  
Therefore, $d_{i}d_{i+1} = 0$ for every $i \in \mathbb{Z}^{+}$, and thus the projective resolution P and the deleted projective resolution $P_{A}$ are both complexes.  
However, $P_{A}$ is not necessarily exact since $im({d_{1}}) = \ker{(\epsilon)}$, which may not equal the kernel of the morphism $P_{0} \rightarrow 0$.  
Now, we can form the {\it induced complex} $TP_{A}$, which is defined as $\cdot\cdot\cdot \rightarrow T(P_{2}) \xrightarrow{T(d_{2})} T(P_{1}) \xrightarrow{T(d_{1})} T(P_{0}) \rightarrow 0$.

For $n \in \mathbb{Z}$, the {\it $n^{th}$ homology} is $H_{n}(C) = Z_{n}(C)/B_{n}(C)$, where $C$ is a complex, $Z_{n}(C) = \ker{(d_{n})}$, and $B_{n}(C) = im(d_{n+1})$.  
Hence, $H_{n}(C) = \ker{(d_{n})}/im(d_{n+1})$.  
If we consider the deleted projective resolution $P_{A}$ as defined above, then $\cdot\cdot\cdot \rightarrow P_{2} \bigotimes_{R} B \xrightarrow{d_{2} \otimes 1_{B}} P_{1} \bigotimes_{R} B \xrightarrow{d_{1} \otimes 1_{B}} P_{0} \bigotimes B \rightarrow 0$ is the induced complex $T_{B}P_{A}$ of the functor $T_{B}(\underline{\hspace{.3cm}})= \underline{\hspace{.3cm}} \bigotimes_{R} B$.
Define the {\it Tor functor} to be $Tor_{n}^{R}(A, B) = H_{n}(T_{B}P_{A}) = \ker{(d_{n} \otimes 1_{B})}/im(d_{n+1} \otimes 1_{B})$.
Note that $Tor_{n}^{R}(A, B)$ does not depend on the choice of projective resolution \cite{Rotman}.
The functor $Tor_{n}^{R}(A,\underline{\hspace{.3cm}})$ is referred to as the {\it left derived functor} of $A \bigotimes_{R} B$.
The following two well-known propositions will be useful later:

\begin{prop}\cite{Rotman}
If $M \in Mod_{R}$ and $0 \rightarrow A \rightarrow B \rightarrow C \rightarrow 0$ is an exact sequence of left $R$-modules, then the induced sequence $...\rightarrow Tor_{n+1}^{R}(M,C) \rightarrow Tor_{n}^{R}(M,A) \rightarrow Tor_{n}^{R}(M,B) \rightarrow Tor_{n}^{R}(M,C) \rightarrow ... \rightarrow Tor_{1}^{R}(M,C) \rightarrow M \bigotimes_{R} A \rightarrow M \bigotimes_{R} B \rightarrow M \bigotimes_{R} C \rightarrow 0$ is exact.
\end{prop}

\begin{prop}\cite{Rotman}
A right $R$-module $M$ is flat if and only if $Tor_{n}^{R}(M,X)=0$ for every left $R$-module $X$ and every $n \geq 1$.
\end{prop}

\chapter{Torsion-free Rings and Modules}

In 1960, Hattori used the homological properties of classical torsion-free modules over integral domains to give a more general definition of torsion-freeness.  
He defines a right $R$-module $M$ to be {\it torsion-free} if $Tor_{1}^{R}(M, R/Rr) = 0$ for every $r \in R$, and he defines a left $R$-module $N$ to be {\it torsion-free} if $Tor_{1}^{R}(R/sR, N) = 0$ for every $s \in R$ \cite{Hattori}.
The following equivalent definition of torsion-freeness is also given by Hattori in \cite[Proposition 1]{Hattori}:

\begin{prop}\label{Hattori Torsion}\cite{Hattori}
The following are equivalent for a right $R$-module $M$.
\begin{enumerate}[(a)]
\item $M$ is torsion-free 
\item For each $x \in M$ and $r \in R$, $xr = 0$ implies the existence of $x_{1},x_{2},...,x_{k} \in M$ and $r_{1},r_{2},...r_{k} \in R$ such that $x = \displaystyle\sum_{j=1}^{k} x_{j}r_{j}$ and $r_{j}r = 0$ for every $j = 1,2,...,k$.
\end{enumerate}
\begin{proof}
Consider the exact sequence $0 \rightarrow Rr \xrightarrow{\iota} R \xrightarrow{\pi} R \slash Rr \rightarrow 0$ of left $R$-modules, where $\iota$ is the inclusion map and $\pi$ is the epimorphism $r \mapsto r+Rr$.
This induces a long exact sequence $X =...\rightarrow Tor_{1}^{R}(M,R \slash Rr) \xrightarrow{f} M \bigotimes_{R} Rr \xrightarrow{1_{M} \otimes \iota} M \bigotimes_{R} R \cong M \xrightarrow{1_{M} \otimes \pi} M \bigotimes_{R} R \slash Rr \rightarrow 0$ \cite[Corollary 6.30]{Rotman}.
Observe that condition ($b$) is equivalent to $1_{M} \otimes \iota$ being a monomorphism.
For if $1_{M} \otimes \iota:x\otimes r \mapsto xr$ is a monomorphism, then $xr=0$ implies $x \otimes r=0$.
Hence, there exists $x_{1}, x_{2},...,x_{k} \in M$ and $r_{1},r_{2},...,r_{k} \in R$ such that $x=x_{1}r_{1}+x_{2}r_{2}+...+x_{k}r_{k}$ and $r_{j}r=0$ for $j=1,2,...,k$ by \Cref{tensor zero}.
On the other hand, if $xr=0$ implies $x=x_{1}r_{1}+x_{2}r_{2}+...+x_{k}r_{k}$ and $r_{j}r=0$, then $x \otimes r=x_{1}r_{1}+x_{2}r_{2}+...+x_{k}r_{k} \otimes r=x_{1}\otimes r_{1}r+x_{2}\otimes r_{2}r+...+x_{k} \otimes r_{k}r=0$.
Hence, $\ker{(1_{M} \otimes \iota)}=0$ and $1_{M} \otimes \iota$ is a monomorphism.

To complete the proof, it needs to be shown that $M$ is torsion-free if and only if $1_{M} \otimes \iota$ is a monomorphism.
If $M$ is torsion-free, then $Tor_{1}^{R}(M,R \slash Rr)=0$.
Thus, $0 \rightarrow M \bigotimes_{R} Rr \xrightarrow{1_{M} \otimes \iota} M \bigotimes_{R} R \cong M \xrightarrow{1_{M} \otimes \pi} M \bigotimes_{R} R \slash Rr \rightarrow 0$ is exact and so $1_{M} \otimes \iota$ is a monomorphism.
Conversely, if $1_{M} \otimes \iota$ is a monomorphism, then $im(f)=\ker{(1_{M} \otimes \iota)}=0$ in the induced sequence $X$.
However, $f$ is a monomorphism.
Hence, $0=im(f) \cong Tor_{1}^{R}(M,R \slash Rr)$.
\end{proof}
\end{prop}

A ring $R$ is {\it torsion-free} if every finitely generated right (or left) ideal is torsion-free as a right (or left) $R$-module.  
Hattori shows in \cite{Hattori} that a ring R is torsion-free if and only if every principal left ideal of R is flat.  
To see this, observe that if $0 \rightarrow J \xrightarrow{i} R \xrightarrow{p} R/J \rightarrow 0$ is an exact sequence of right R-modules with $J$ finitely generated, then $0 \rightarrow J \bigotimes_{R} Rr \xrightarrow{i \otimes 1_{Rr}} R \bigotimes_{R} Rr \xrightarrow{p \otimes 1_{Rr}} R/J \bigotimes_{R} Rr \rightarrow 0$ is an exact sequence whenever $Rr$ is flat.
This is the case if and only if $Tor_{1}^{R}(R/J, Rr) = 0$.
Hattori gives a natural isomorphism in \cite[Proposition 7]{Hattori} showing that $Tor_{1}^{R}(R/J, Rr) \cong Tor_{1}^{R}(J,R \slash Rr)$.
Hence, $Tor_{1}^{R}(J,R \slash Rr)=0$ if and only if $Rr$ is flat for every $r \in R$.
That is, every finitely generated right ideal is torsion-free if and only if every principal left ideal is flat.

In 2004, John Dauns and Lazlo Fuchs provided the following useful characterization of torsion-free rings:

\begin{theorem}\label{Dauns Fuchs 3.2}\cite{Dauns}The following are equivalent for a ring R:
\begin{enumerate}[(a)]
   \item R is torsion-free.
   \item For every $s, r \in R$, $sr = 0$ if and only if $s \in s \cdot ann_{l}(r)$.  In other words, $sr = 0$ if and only if $s = su$ and $ur = 0$ for some $u \in R$.      
  \end{enumerate}

\begin{proof}(a) $\Rightarrow$ (b): Suppose R is a torsion-free ring.  
For $s \in R$, $sR$ is torsion-free as a right R-module.  
By \Cref{Hattori Torsion}, if $a \in sR$ and $r \in R$ with $ar = 0$, then there exists $u \in R$ so that $a = su$ and $ur = 0$.
Hence, if $sr = 0$,  we have $s = su$ and $ur = 0$ for some $u \in R$, since $s = s\cdot1$ $\in sR$.  
Conversely, if there is some $u \in R$ such that $s = su$ and $ur = 0$, then $sr = (su)r = s(ur) = s\cdot0 = 0$.  
Therefore, $sr = 0$ if and only if $s = su$ and $ur = 0$ for some $u \in R$.

(b) $\Rightarrow$ (a): Assume that $sr = 0$ for every $s, r \in R$ if and only if $s = su$ and $ur = 0$ for some $u \in R$.
Let $Rr$ be a finitely generated left ideal of R.  
Assume that the sequence $0 \rightarrow J \rightarrow R \rightarrow R/J \rightarrow 0$ is exact with $J$ finitely generated.
Then, R is a torsion-free ring if $0 \rightarrow J \bigotimes_{R} Rr \xrightarrow{\varphi} R \bigotimes_{R} Rr \xrightarrow{\psi} R/J \bigotimes_{R} Rr \rightarrow 0$ is exact.
By \Cref{Right Exactness}, it follows that $J \bigotimes_{R} Rr \xrightarrow{\varphi} R \bigotimes_{R} Rr \xrightarrow{\psi} R/J \bigotimes_{R} Rr \rightarrow 0$ is exact.
In order for the entire sequence to be exact, it needs to be shown that $\varphi$ is a monomorphism.

Note that $R \bigotimes_{R} Rr \cong Rr$ by \Cref{Tensor Nat Iso}.
Consider $j \otimes sr \in J \bigotimes_{R} Rr$.  
Since $j \otimes sr = js \otimes r$ and $js \in J$, tensors in $J \bigotimes_{R} Rr$ can be written as $k \otimes r$ for some $k \in J$.
Thus, it needs to be shown that $J \bigotimes_{R} Rr \xrightarrow{\varphi} Rr$ given by $\varphi(k \otimes r) = kr$ is a monomorphism.
Let $k \otimes r \in \ker{\varphi}$.  Then $\varphi(k \otimes r) = kr = 0$.
By assumption, there exists some $u \in R$ such that $k = ku$ and $ur = 0$.
Then, $k \otimes r = ku \otimes r = k \otimes ur = k \otimes 0 = 0$.
Thus, $\ker{\varphi} = 0$ and $\varphi$ is a monomoprhism.
Therefore, $0 \rightarrow J \bigotimes_{R} Rr \xrightarrow{\varphi} R \bigotimes_{R} Rr \xrightarrow{\psi} R/J \bigotimes_{R} Rr \rightarrow 0$ is an exact sequence, and hence $R$ is a torsion-free ring.
\end{proof}
\end{theorem}

\begin{prop}\cite[Proposition 7]{Hattori}\label{submod of tor free ring}
A ring $R$ is torsion-free if and only if every submodule of a torsion-free right $R$-module is torsion-free.
\begin{proof}
Suppose $R$ is torsion-free and let $N$ be a submodule of a torsion-free right $R$-module $M$.
Consider the exact sequence $0 \rightarrow N \xrightarrow{\iota} M \xrightarrow{\pi} M \slash N \rightarrow 0$, where $\iota$ is the inclusion map and $\pi$ is the canonical epimorphism.
As noted above, if $R$ is torsion-free, then the principal left ideal $Rr$ is flat for every $r \in R$.
Hence, $0 \rightarrow N \bigotimes_{R} Rr \rightarrow M \bigotimes_{R} Rr \rightarrow M \slash N \bigotimes_{R} Rr \rightarrow 0$ is exact and so $Tor_{1}^{R}(M\slash N,Rr)\cong 0$.
Observe that $Tor_{1}^{R}(M,R\slash Rr) \cong 0$ since $M$ is torsion-free.
If we consider the long exact sequence derived from the functor $Tor_{n}^{R}(\underline{\hspace{.3cm}},R \slash Rr)$, then $0\cong Tor_{1}^{R}(M\slash N,Rr) \cong Tor_{2}^{R}(M\slash N,R \slash Rr) \rightarrow Tor_{1}^{R}(N,R \slash Rr) \rightarrow Tor_{1}^{R}(M,R \slash Rr) \cong 0$ is exact.
Therefore, $Tor_{1}^{R}(N,R \slash Rr)=0$ and $N$ is torsion-free.
On the other hand, if every submodule of a torsion-free right $R$-module is torsion-free, then every finitely generated right ideal of $R$ is torsion-free since $R$ itself is torsion-free as a right $R$-module. 
\end{proof}
\end{prop}

\begin{theorem}\label{Dauns Fuchs 4.5}\cite{Dauns}A ring $R$ is a right p.p.-ring if and only if $R$ is torsion-free and, for each $x \in R$, $ann_{r}(x)$ is finitely generated.
\begin{proof}Suppose $R$ is a right p.p.-ring.  
Then, for each $r \in R$, $ann_{r}(r) = eR$ for some idempotent $e \in R$.
Let $s \in R$ be such that $rs = 0$.  Then, $s \in ann_{r}(r)$, and hence $s = es^{'}$ for some $s^{'} \in R$.
It follows that $es = e^{2}s^{'} = es^{'} = s$.
Furthermore, $e = e^{2} \in eR = ann_{r}(r)$ and hence $re = 0$.
Note also that if $s = es$ and $re = 0$, then $s \in eR = ann_{r}(r)$ and hence $rs = 0$.
Thus, $rs = 0$ if and only if $s = es$ and $re = 0$.
Therefore, $R$ is a torsion-free ring by a symmetric version of \Cref{Dauns Fuchs 3.2}.
Moreover, since $R$ is a right p.p.-ring, $ann_{r}(r)$ is generated by an idempotent and thus finitely generated.

Conversely, suppose $R$ is a torsion-free ring and the right annihilator of every element of $R$ is finitely generated.
Let $s \in R$ and let $\{s_{1},...,s_{n}\}$ be the finite set of generators for $ann_{r}(s)$.
Note that each $s_{i} \in ann_{r}(s)$, and so $ss_{i} = 0$ for each $i = 1,...,n$.
Let $S = \displaystyle \bigoplus^{n} R$ be the direct sum of n copies of $R$, and consider $S$ as a left $R$-module.
Let $s^{'} = (s_{1},...s_{n}) \in S$.
Note that $S$ is a torsion-free left $R$-module since it is the direct sum of copies of $R$, which is torsion-free as a left $R$-module.
Thus, the submodule $Rs^{'}$ of $S$ is torsion-free by \Cref{submod of tor free ring}.
Hence, \Cref{Hattori Torsion} gives some $u \in R$ such that $s^{'} = us^{'}$ and $su = 0$, and thus $u \in ann_{r}(s)$.  
Note that $s_{i} = us_{i}$ for each $i = 1,...,n$.
This implies that $s_{i} \in uR$ for each $i$, and so $\{s_{1},...,s_{n}\} \subseteq uR$.
It follows that $ann_{r}(s) = s_{1}R + ... + s_{n}R \subseteq uR$.
Suppose $x \in uR$.  Then, $x = ut$ for some $t \in R$.
Thus, $sx = sut = 0\cdot t = 0$, and so $x \in ann_{r}(s)$.
Therefore, $ann_{r}(s) = uR$.

Now, since $R$ is a torsion-free ring, $uR$ is torsion-free as a finitely generated right ideal of $R$.
By a symmetric version of \Cref{Dauns Fuchs 3.2}, since $su = 0$, there exists an $e \in uR = ann_{r}(s)$ such that $u = eu$ and $se = 0$.
Let $x \in uR$.  
Then, $x = ut = eut \in eR$ for some $t \in R$.  
Hence, $uR \subseteq eR$.
On the other hand, suppose $y \in eR$.  Then, for some $v \in R$, $y = ev$ and $sy = sev = 0\cdot v = 0$.
Thus, $y \in ann_{r}(s)$ and $eR \subseteq ann_{r}(s) = uR$.
Hence, $ann_{r}(s) = uR = eR$ and $e = ur$ for some $r \in R$.
It then follows that $e$ is an idempotent since $e^{2} = eux = ux = e$.
Therefore, $ann_{r}(s)$ is generated by an idempotent and so $R$ is a right p.p.-ring.
\end{proof}
\end{theorem}

\begin{lemma}\label{eR is ann}
If $R$ is a right p.p.-ring and $e \in R$ is a nonzero idempotent, then $eR = ann_{r}(x)$ for some $x \in R$.  In particular, $eR=ann_{r}(1-e)$.
\begin{proof}
If $er \in eR$, then $(1-e)er = (e - e^2)r = (e - e)r = 0$.
Hence, $er \in ann_{r}(1-e)$ and $eR \subseteq ann_{r}(1-e)$.
On the other hand, if $s \in ann_{r}(1-e)$, then $(1-e)s = 0$.
Hence, $s - es = 0$, and so $s = es \in eR$.
Therefore, $eR = ann_{r}(1-e)$.
\end{proof}
\end{lemma}

\begin{prop}\cite{Albrecht}\label{ADF 3.3}
If $R$ is a right and left p.p.-ring which does not contain an infinite set of orthogonal idempotents and $M$ is a torsion-free right R-module, then $ann_{r}(x)$ is generated by an idempotent for every $x \in M$.
\begin{proof}
Let $R$ be a right and left p.p.-ring which does not contain an infinite set of orthogonal idempotents.
Take $M$ to be a torsion-free right $R$-module and let $A=ann_{r}(x)$ for some nonzero $x \in M$.
Suppose $r_{0} \in R$ is such that $xr_{0} = 0$.
Note that the cyclic submodule $xR$ is torsion-free since $R$ is a right p.p.-ring.
Moreover, $ann_{l}(r_{0})= Re_{0}$ for some idempotent $e_{0} \in R$ since $R$ is a left p.p.-ring.
By \Cref{Hattori Torsion}, there exists $xs_{1},xs_{2},...,xs_{n} \in xR$ and $t_{1}e_{0},t_{2}e_{0}...,t_{n}e_{0} \in Re_{0}=ann_{l}(r_{0})$ such that $x = xs_{1}t_{1}e_{0}+xs_{2}t_{2}e_{0}+...+xs_{n}t_{n}e_{0}$.
Hence, $xe_{0}=xs_{1}t_{1}e_{0}^2+xs_{2}t_{2}e_{0}^2+...+xs_{n}t_{n}e_{0}^2 = x$.
Thus, $0=x-xe_{0}=x(1-e_{0})$.
Therefore, if $(1-e_{0})r \in (1-e_{0})R$, then $x(1-e_{0})r =0$ and $(1-e_{0})R \subseteq A$.

Now, if there exists some $r_{1} \in A \backslash (1-e_{0})R$, then $r_{1} \neq (1-e_{0})r_{1}$ and hence $e_{0}r_{1} \neq 0$.
However, $xe_{0}r_{1}=xr_{1}=0$.
Since $R$ is a left p.p.-ring, $ann_{l}(e_{0}r_{1})=R(1-f)$ for some idempotent $1-f$.
Note that as before it follows from \Cref{Hattori Torsion} that $x=x(1-f)$ since $xe_{0}r_{1}=0$.
Furthermore, $1-e_{0} \in ann_{l}(e_{0}r_{1})=R(1-f)$ since $(1-e_{0})e_{0}r_{1}=e_{0}r_{1}-e_{0}r_{1}=0$.
Hence, there is some $r \in R$ such that $(1-e_{0})f=r(1-f)f=r(f-f)=0$.
Thus, $e_{0}f=f$. 
Let $e_{1}=(1-f)e_{0}=e_{0}-fe_{0}$.
Then, $e_{1}^2 = (e_{0}-fe_{0})(e_{0}-fe_{0})=e_{0}-e_{0}fe_{0}-fe_{0}+fe_{0}fe_{0}=e_{0}-fe_{0}-fe_{0}+fe_{0}=e_{0}-fe_{0} = e_{1}$.
Thus, $e_{1}$ is an idempotent.
Moreover, $e_{1}$ is nonzero, since otherwise $e_{0}=fe_{0}$ and hence $e_{0}=0$.

Now, $e_{1}e_{0}=(1-f)e_{0}e_{0}=(1-f)e_{0}=e_{1}$, and \Cref{eR is ann} shows that $(1-e_{0})R = ann_{r}(e_{0})$ and $(1-e_{1})R = ann_{r}(e_{1})$.
Thus, if $r \in ann_{r}(e_{0})$, then $e_{1}r=e_{1}e_{0}r=0$.
Hence, $r \in ann_{r}(e_{1})=(1-e_{1})R$, and so $(1-e_{0})R \subseteq (1-e_{1})R$.
Moreover, $e_{1}e_{0}r_{1}=e_{1}r_{1}=(1-f)e_{0}r_{1}=0$ since $1-f \in ann_{r}(e_{0}r_{1})$.
Thus, $e_{0}r_{1} \in ann_{r}(e_{1})=(1-e_{1})R$.
However, $e_{0}r_{1}$ is nonzero and hence $e_{0}r_{1} \notin ann_{r}(e_{0})=(1-e_{0})R$.
Thus, $(1-e_{0})R \subset (1-e_{1})R$ is a proper inclusion.
By supposing there is some $r_{2} \in A \backslash (1-e_{1})R$ and repeating these steps, and then supposing there is some $r_{3} \in A \backslash (1-e_{2})R$ and so on, we can construct an ascending chain $(1-e_{0})R \subset (1-e_{1})R \subset (1-e_{2})R \subset ...$.
However, this chain must terminate at some point since $R$ only contains finite sets of orthogonal idempotents.  
Therefore, there is some idempotent $e \in R$ such that $A = (1-e)R$.
\end{proof}
\end{prop}

\begin{prop}\cite{Albrecht}\label{AFD 3.4}
If $R$ is a right and left p.p.-ring not containing an infinite set of orthogonal idempotents, then a cyclic submodule of a torsion-free right $R$-module is projective.
\begin{proof}
Let $M$ be a torsion-free right $R$-module, and take $N$ to be a cyclic submodule of $M$.
Then, $N$ is of the form $xR$ for some $x \in N \leq M$.
By \Cref{ADF 3.3}, $ann_{r}(x)=eR$ for some idempotent $e \in R$.
If $f:R \rightarrow xR$ is the epimorphism defined by $r \mapsto xr$, then $xR \cong R \slash \ker{(f)} = R \slash ann_{r}(x)$ by the First Isomorphism Theorem.
It then follows that $xR \cong R \slash ann_{r}(x) \cong [eR \bigoplus (1-eR)] \slash ann_{r}(x) \cong$ $[eR \bigoplus (1-e)R] \slash eR \cong (1-e)R$.
Therefore, $N$ is a principal right ideal of $R$, and thus projective, since $R$ is a right p.p.-ring. 
\end{proof} 
\end{prop}

A ring $R$ is a {\it Baer-ring} if $ann_{r}(A)$ is generated by an idempotent for every subset $A$ of $R$.
Note that if $R$ is Baer, then $ann_{r}(ann_{l}(A))=eR$ for some idempotent $e \in R$.
Hence, $ann_{l}(A)=ann_{l}(ann_{r}(ann_{l}(A)))=ann_{l}(eR)=R(1-e)$ by \Cref{eR is ann}.
Thus, $ann_{r}(A)$ is generated by an idempotent if and only if $ann_{l}(A)$ is generated by an idempotent.
Therefore, the property that $R$ is a Baer ring is right-left-symmetric.
The following theorem from Dauns and Fuchs \cite{Dauns} gives conditions for which a ring $R$ is Baer:

\begin{theorem}\cite{Dauns}\label{Fuchs and Dauns 4.8} If R is a torsion-free ring and right annihilators of elements are finitely generated and satisfy the ascending chain condition, then R is a Baer-ring.
\begin{proof}It follows from $\Cref{Dauns Fuchs 4.5}$ that R is a right p.p.-ring since $ann_{r}(x)$ is finitely generated for every $x \in R$. 
Thus, for each $x \in R$, there is some idempotent $e \in R$ such that $ann_{r}(x) = eR$. 
Suppose R contains an infinite set $E$ of orthogonal idempotents.
Consider two idempotents $e_{1}$ and $e_{2}$ in $E$, and let $e_{1}r \in e_{1}R$.
Note that since $e_{1}$ and $e_{2}$ are orthogonal idempotents, $e_{1}r = (e_{1}+0)r = (e_{1}^2 + e_{2}e_{1})r = (e_{1}+e_{2})e_{1}r \in (e_{1}+e{2})R$.
Therefore, $e_{1}R \subseteq (e_{1}+e_{2})R$.
Inductively, we can construct an ascending chain of principal ideals generated by idempotents. 
For if $e_{1}, ..., e_{n}, e_{n+1}$ are orthogonal idempotents in the infinite set and $(e_{1}+...+e_{n})r \in (e_{1}+...+e_{n})R $, then 
$(e_{1}+e_{2}+...+e_{n})r = (e_{1}^2+e_{2}^2...+e_{n}^2+0)r =[(e_{1}^2+e_{1}e_{2}+...e_{1}e_{n})+(e_{2}e_{1}+e_{2}^2+...+e_{2}e_{n})+...+(e_{n}e_{1}+...+e_{n}^2)+(e_{n+1}e_{1}+...+e_{n+1}e_{n})]r = (e_{1}+...+e_{n+1})(e_{1}+...+e_{n})r \in (e_{1}+...+e_{n+1})R$.

Hence, $e_{1}R \subseteq (e_{1}+e_{2})R \subseteq ... \subseteq (e_{1}+...+e_{n})R \subseteq (e_{1}+...+e_{n+1})R \subseteq ...$ 
is an ascending chain of principal ideals generated by idempotents.
Furthermore, this will be an infinite chain since there are an infinite number of idempotents in $E$.
Note that by \Cref{eR is ann}, for each $n \in \mathbb{Z}^{+}$, $(e_{1}+...+e_{n})R = ann_{r}(x)$ for some $x \in R$.
Thus, an infinite ascending chain of right annihilators has been constructed, contradicting the ascending chain condition on right annihilators.
Therefore, $R$ does not contain an infinite set of orthogonal idempotents.
Since R is a right p.p.-ring which does not contain an infinite set of orthogonal idempotents, by \Cref{RwCC 8.4} every right annihilator in R is generated by an idempotent.
Therefore, R is a Baer-ring.
\end{proof}
\end{theorem}

\chapter{Non-singularity}

\section{Essential Submodules and the Singular Submodule}

Let $R$ be a ring and consider a submodule $A$ of a right $R$-module $M$.  
If $A \cap B$ is nonzero for every nonzero submodule $B$ of $M$, then $A$ is said to be an {\it essential} submodule of $M$.  This is denoted $A \leq^{e} M$.  
In other words, $A \leq^{e} M$ if and only if $B = 0$ whenever $B \leq M$ is such that $A \cap B = 0$.
A monomorphism $\alpha:A \rightarrow B$ is called {\it essential} if $im (A) \leq^{e} B$.

\begin{prop}\cite[Corollary 5.13]{Anderson}\label{And 5.13}
A monomorphism $\alpha:A \rightarrow B$ is essential if and only if, for every right $R$-module $C$ and every $\beta \in Hom_{R}(B,C)$, $\beta$ is a monomorphism whenever $\beta\alpha$ is a monomorphism.
\end{prop}

The {\it singular submodule of $M$} is defined as $\textbf{Z}(M) = \{x \in M$ $|$ $xI = 0$ for some essential right ideal $I$ of R\}.
Equivalently, $\textbf{Z}(M) = \{x \in M$ $|$ $ann_{r}(x) \leq^{e} R$\}.
For if $I \leq^{e} R$ and $x \in M$ is such that $xI = 0$, then for any nonzero right ideal $J$ of $R$, there is an element $a \in I \cap J$.  
Since $a \in I$,  $xa = 0$. Hence, $a \in ann_{r}(x) \cap J$ and so $ann_{r}(x) \leq^{e} R$. 
On the other hand, note that $ann_{r}(x)$ is a right ideal of $R$ such that $x \cdot ann_{r}(x) = 0$. 
A right $R$-module $M$ is called {\it singular} if $\textbf{Z}(M) = M$ and {\it non-singular} if $\textbf{Z}(M) = 0$.
If $R$ is viewed a right $R$-module, then the {\it right singular ideal of $R$} is $\textbf{Z}_{r}(R) = \textbf{Z}(R_{R})$.
The ring $R$ is {\it right non-singular} if it is non-singular as a right $R$-module.

\begin{prop}\cite{Goodearl}\label{Goodearl 1.20a}
A right $R$-module $A$ is non-singular if and only if $Hom_{R}(C,A) = 0$ for every singular right $R$-module $C$.
\begin{proof}
Suppose $A$ is a non-singular right $R$-module and $C$ is a singular right $R$-module.
Let $f \in Hom_{R}(C,A)$.
If it can be shown that $f(\textbf{Z}(C)) \leq \textbf{Z}(A)$, then the proof follows readily since $f(C) = f(\textbf{Z}(C))$ and $\textbf{Z}(A) = 0$.
Suppose $x \in \textbf{Z}(C)$.
Then, $ann_{r}(x) \leq^{e} R$.
Hence, if $I$ is any nonzero right ideal of $R$, then there exists some $y \in I$ such that $xy = 0$.
Then, $f(x)y = f(xy) = f(0) = 0$ and $y \in ann_{r}(f(x)) \cap I$.
Thus, $ann_{r}(f(x)) \leq^{e} R$ and so $f(x) \in \textbf{Z}(A)$.
Therefore, $f(\textbf{Z}(C)) \leq \textbf{Z}(A)$.

Conversely, suppose $A$ is a right $R$-module and $Hom_{R}(C,A) = 0$ for every singular right $R$-module $C$.
Then, $Hom_{R}(\textbf{Z}(A),A) = 0$ since the singular submodule $\textbf{Z}(A)$ is singular.
Hence, the inclusion map $\iota: \textbf{Z}(A) \rightarrow A$ given by $\iota(x) = x$ is a zero map.
Thus, $\textbf{Z}(A) = \iota(\textbf{Z}(A)) = 0$.
Therefore, $A$ is a non-singular right $R$-module. 
\end{proof}
\end{prop}

\begin{prop}\cite{Goodearl}\label{Goodearl 1.20b}
The following are equivalent for a right $R$-module $C$:
\begin{enumerate}[(a)]
\item $C$ is singular.
\item There exists an exact sequence $0 \rightarrow A \xrightarrow{f} B \xrightarrow{g} C \rightarrow 0$ such that $f$ is essential. 
\end{enumerate}
\begin{proof}
($a$) $\Rightarrow$ ($b$): Suppose $C$ is a right $R$-module.
Let $0 \rightarrow A \xrightarrow{\iota} B \xrightarrow{g} C \rightarrow 0$ be an exact sequence of right $R$-modules such that $B$ is free and $\iota$ is the inclusion map. 
Let $\{x_{\alpha}\}_{\alpha \in K}$ be a basis for $B$ for some index $K$.
Then, for each $\alpha \in K$, $g(x_{\alpha}) \in C = \textbf{Z}(C)$.
Hence, there exists an essential right ideal $I_{\alpha}$ of $R$ such that $g(x_{\alpha}I_{\alpha}) = g(x_{\alpha})I_{\alpha} = 0$.
Thus, for each $\alpha \in K$ and each $i_{\alpha} \in I_{\alpha}$, $x_{\alpha}i_{\alpha} \in \ker{g} = A$.
That is, $x_{\alpha}I_{\alpha} \leq A$ for each $\alpha \in K$, and it follows that $\bigoplus_{K} x_{\alpha}I_{\alpha} \leq A$.
If $x_{\alpha}J$ is a nonzero right ideal of $x_{\alpha}R$, then $J$ is a nonzero right ideal of $R$, and there is a nonzero element $y \in I_{\alpha} \cap J$.
Then it readily follows that $x_{\alpha}y \in x_{\alpha}I_{\alpha} \cap x_{\alpha}J$ is nonzero.
Hence, $x_{\alpha}I_{\alpha} \leq^{e} x_{\alpha}R$ for each $\alpha \in K$. 
Thus, $\bigoplus_{K} x_{\alpha}I_{\alpha} \leq^{e} \bigoplus_{K} x_{\alpha}R = B$.
Therefore, $A$ is also essential in $B$ since $\bigoplus_{K} x_{\alpha}I_{\alpha} \leq A$.  
It then follows from the exactness of the sequence that $im(A) \cong A \leq^{e} B$.

($b$) $\Rightarrow$ ($a$): Assume $0 \rightarrow A \xrightarrow{f} B \xrightarrow{g} C \rightarrow 0$ is an exact sequence of right $R$-modules such that $im(A) \leq^{e} B$.
For each $b \in B$, define $h_{b}: R \rightarrow B$ by $h_{b}(r) = br$, and let \\$I_{b} = \{r \in R$ $|$ $br \in im(A)\}$.
Note that $I_{b}$ is a nonzero right ideal of $R$.
Suppose $I_{b}$ is not essential in $R$.
Then there is a nonzero right ideal $J$ of $R$ such that $I_{b} \cap J = 0$.
Moreover, if $s \in \ker{(h_{b})}$, then $h_{b}(s) = bs = 0 \in im(A)$ and it follows that $\ker{(h_{b})} \subseteq I_{b}$. 
Hence, $\ker{(h_{b})} \cap J = 0$.
Thus, $h_{b}|_{J}$ is a monomorphism.
This implies that $h_{b}(J)$ must be a nonzero right ideal of $B$ since $J$ is a nonzero right ideal of $R$.
Thus, $h_{b}(J) \cap im(A) \neq 0$ by the assumption that $im(A) \leq^{e} B$.
Then for some nonzero $j \in J$, $bj = h_{b}(j) \in im(A)$.
Hence, $j \in I_{b} \cap J$, which is a contradiction. 
Therefore, $I_{b}$ is an essential right ideal of $R$.
Note that for every $b \in B$, if $bi \in bI_{b}$, then $bi \in im(A)$.
Then by exactness of the sequence, $bI_{b} \subseteq im(A) = \ker{g}$.
Hence, $g(b)I_{b} = g(bI_{b}) = 0$, which implies $g(b) \in \textbf{Z}(C)$.
Since this is the case for every $b \in B$, $g(B) \subseteq \textbf{Z}(C)$.
Furthermore, since the sequence is exact, $C = g(B) \subseteq \textbf{Z}(C)$.
Therefore, $C = \textbf{Z}(C)$.
\end{proof}
\end{prop}

\begin{prop}\label{pp implies ns}
If R is a right p.p.-ring, then R is a right non-singular ring.
\begin{proof}
Let $R$ be a right p.p.-ring and take any $x \in R$.
Suppose $ann_{r}(x) \leq^{e} R$.
Since $R$ is a right p.p.-ring, $ann_{r}(x)=eR$ for some idempotent $e \in R$.
Observe that $R= eR \bigoplus (1-e)R$.
Hence, $ann_{r}(x) \cap (1-e)R=0$.
However, this implies that $(1-e)R = 0$ since $ann_{r}(x) \leq^{e} R$.
Hence, $1 - e = 0$, and so $ann_{r}(x) = 1R = R$.
Thus, $xr = 0$ for every $r \in R$, which implies $x = 0$.
Therefore, $R$ is right non-singular.
\end{proof}
\end{prop}

\section{The Maximal Ring of Quotients and Right Strongly Non-singular Rings}

The maximal ring of quotients and strongly non-singular rings will play an important role in determining which rings satisfy the condition that the classes of torsion-free and non-singular modules coincides.
We explore these concepts in this section.
If $R$ is a subring of a ring $Q$, then $Q$ is a {\it classical right ring of quotients} of $R$ if every regular element of $R$ is a unit in $Q$ and every element of $Q$ is of the form $rs^{-1}$, where $r,s \in R$ with $s$ regular \cite{Hungerford}.
For a ring which is not necessarily commutative, such a $Q$ may not exist.
Thus, we consider a more general way to define the right ring of quotients which guarantees its existence for any ring $R$.

Let $A$ be a submodule of a right $R$-module $B$.  If $Hom_{R}(M \slash A, B)=0$ for every right $R$-module $M$ satisfying $A \leq M \leq B$, then $B$ is a {\it rational extension} of $A$.
This is denoted $A \leq^{r} B$.

\begin{lemma}\cite{Goodearl}\label{Goodearl 2.24}
Let $B$ be a non-singular right $R$-module and take any submodule $A$ of $B$.
Then, $A \leq^{r} B$ if and only if $A \leq^{e} B$.
\begin{proof}
Suppose $A \leq^{r} B$ and let $M \leq B$ be such that $M \cap A = 0$.
Now, $M \bigoplus A$ is a right $R$-module satisfying $A \leq M \bigoplus A \leq B$.
Hence, $Hom_{R}([M \bigoplus A] \slash A,B)=0$.
Consider $f:(M \bigoplus A) \slash A \rightarrow M$ defined by $(m+a)+A \mapsto m$ for $m \in M$ and $a \in A$.
If $m,m_{0} \in M$ and $a, a_{0} \in A$ are such that $(m+a)+A=(m_{0})+a_{0})+A$, then $(m-m_{0})+(a-a_{0}) \in A$.
Hence, $m-m_{0} \in A$.
However, $M \cap A = 0$ and so $m-m_{0}=0$.
Thus, $f$ is well-defined.
Moreover, $f$ is an isomorphism. 
For if $m \in M$, then $f[(m+a)+A]=m$ for any $a \in A$, and $f[(m+a)+A]=0$ implies that $(m+a)+A=m+A=0$.
Observe that $f \in Hom_{R}([M \bigoplus A] \slash A,B)=0$ since $M \leq B$.
Thus, $M=im(f)=0$ and therefore $A \leq^{e} B$.
Note that this implication does not require $B$ to be right non-singular.

On the other hand, suppose $A \leq^{e} B$ and take $M$ to be a right $R$-module such that $A \leq M \leq B$.
Then, any nonzero submodule $N$ of $B$ is such that $A \cap N \neq 0$.
Hence, any nonzero submodule $K$ of $M$ is such that $A \cap K \neq 0$ since any such submodule is also a submodule of $B$.
Thus, $A \leq^{e} M$.
Consider the exact sequence $0 \rightarrow A \xrightarrow{\iota} M \xrightarrow{\pi} M \slash A \rightarrow 0$, where $\iota$ is the inclusion map and $\pi$ is the canonical epimorphism.
Observe that $im(\iota)=A \leq^{e} M$.
Hence, $M/A$ is singular by \Cref{Goodearl 1.20b}.
It then follows from \Cref{Goodearl 1.20a} that $Hom_{R}(M \slash A,B)=0$ since $B$ is nonsingular.
Therefore, $B$ is a rational extension of $A$.
\end{proof}
\end{lemma}

A right $R$-module $E$ is called {\it injective} if, given any two right $R$-modules $A$ and $B$, a monomorphism $\alpha:A \rightarrow B$, and a homomorphism $\varphi:A \rightarrow E$, there exists a homomorphism $\psi:B \rightarrow E$ such that $\varphi = \psi\alpha$.
If $E$ is injective and $M_{R} \leq^{e} E_{R}$, then $E$ is called an {\it injective hull} of $M$.
Every right $R$-module $M$ has an injective hull, which is unique up to isomorphism \cite[Theorems 1.10, 1.11]{Goodearl}.

Let $R$ be a subring of a ring $Q$.
If $R_{R} \leq^{r} Q_{R}$, then $Q$ is a {\it right ring of quotients of $R$}. 
Observe that $R$ is a right ring of quotients of itself since $R_{R} \leq^{r} R_{R}$.
Similarly, if $_{R}R \leq^{r}$ $_{R}Q$, then $Q$ is a {\it left ring of quotients of $R$}.
Let $Q$ be a right ring of quotients of $R$ such that given any other right ring of quotients $P$ of $R$, the inclusion map $\mu: R \rightarrow Q$ extends to a monomorphism $\nu:P \rightarrow Q$.
Here, $Q$ is called a {\it maximal right ring of quotients of $R$}.
This is denoted $Q^r$ when there is no confusion as to which ring the maximal quotient ring applies, and $Q^r(R)$ otherwise.
The {\it maximal left ring of quotients} $Q^l$ is similarly defined.
In general, $Q^r \neq Q^l$.

\begin{theorem}\cite{Goodearl}\label{Goodearl 2.29}
For any ring $R$, the maximal right ring of quotients $Q^r(R)$ exists.  
In particular, if $E$ is the injective hull of $R_{R}$ and $T=End_{R}(E)$, then $Q = \cap\{\ker{\delta}$ $|$ $\delta \in T$ and $\delta R=0\}$ is a maximal right ring of quotients.
\begin{proof}
If $E$ is the injective hull of $R$, then $\tau x = \tau(x)$ defines a left $T$-module structure on $E$ for $\tau \in T$ and $x \in E$.
Let $T_{0}=End_{T}(E)$ and define $\omega(x)=x\omega$ for $\omega \in T_{0}$ and $x \in E$.
Consider the homomorphisms $\psi:T \rightarrow E$ and $\varphi:T_{0} \rightarrow E$ defined by $\psi\tau = \tau 1$ and $\varphi\omega = 1\omega$.
It is easily seen that $\psi$ is an epimorphism and $\varphi$ is a monomorphism.
Let $x \in E$ and consider the homomorphism $\sigma:R \rightarrow xR$ defined by $\sigma(r) = xr$.
Since $R$ is a subring of $E$, $\sigma$ can be extended to a homomorphism $\tau:E \rightarrow E$.
Thus, $\tau(1)=\sigma(1)=x$ and so $\psi(\tau)=\tau(1)=x$.
Therefore, $\psi$ is an epimorphism.
Now, suppose $\omega \in \ker{\varphi}$.
Then $1\omega=\varphi(\omega)=0$.
If $x \in E$, then $\tau 1=x$ for some $\tau \in T$ since $\psi$ is an epimorphism.
Hence, $\omega(x)=x\omega=(\tau 1)\omega=\tau(1\omega)=\tau(0)=0$.
Therefore, $\omega= 0$ and $\varphi$ is a monomorphism.

If $\delta \in T$ is such that $\delta R=0$, then $\delta(1\omega)=(\delta 1)\omega = 0$ for every $\omega \in T_{0}$.
Hence, $1\omega \in Q$.
Therefore, $\varphi$ can actually be defined as a map $T_{0} \rightarrow Q$.
It readily follows that $\varphi$ maps onto $Q$ and hence $\varphi:T_{0} \rightarrow Q$ is an isomorphism.
To see this, let $x \in Q$ and consider $\nu:E \rightarrow E$ defined by $(\tau 1)\nu = \tau x$. 
This can be defined for every $\tau \in T$ since $\varphi$ is a well-defined epimorphism onto $E$.
Thus, if  $1_{E} \in T$ is the identity map on $E$, then $\varphi(\nu)=1\nu=[1_{E}(1)]\nu=1_{E}(x)=x$.
Therefore, $\varphi$ is onto.

We now define multiplication on $Q$.
For $x, y \in Q$, let $x \cdot y = \varphi[(\varphi^{-1}x)(\varphi^{-1}y)]=1(\varphi^{-1}x)(\varphi^{-1}y)$.
Clearly $x \cdot y \in Q$ and it is easily seen to be associative.
Since $\varphi$ is an isomorphism, if $r \in R$, then there exists some $\omega \in T_{0}$ such that $\varphi(\omega)=1\omega=r$.
Thus, if $x \in Q$, then $x \cdot r = 1(\varphi^{-1}x)(\varphi^{-1}r)=(\varphi\varphi^{-1}x)(\omega)=x\omega=(x1)\omega=x(1\omega)=xr$.
It follows from \cite[Theorem 2.26]{Goodearl} that this multiplication defines a unique ring structure on $Q$ which is consistent with the $R$-module structure.. 

To see that $Q$ is a right ring of quotients, suppose $R \leq M \leq Q$ for some right $R$-module $M$ and let $\alpha \in Hom_{R}(M \slash R,Q)$.
Consider the epimorphism $\pi:M \rightarrow M \slash R$ given by $x \mapsto x+R$, and define $\gamma=\alpha\pi:M \rightarrow Q$.
Observe that $\gamma R=0$ since $\gamma(r)=\alpha\pi(r)=\alpha(r +R)=0$ for any $r \in R$.
Moreover, $\gamma$ can be extended to a map $\beta \in T$ such that $\beta R = 0$.
Since $Q$ is the intersection of the kernels of all homomorphisms $\delta \in T$ satisfying $\delta R = 0$, $M \subseteq Q \subseteq \ker{\beta}$.
Thus, $\gamma M=\beta M=0$ and so $\alpha(x+R)=\gamma(x)=0$ for any $x \in M$.
Therefore, $R \leq^{r} Q$ and $Q$ is a right ring of Quotients.

To see that $Q^r$ is a maximal right ring of quotients, let $P$ be another right ring of quotients. 
Then $R_{R} \leq^{r} P_{R}$ by definition, and hence $R_{R} \leq^{e} P_{R}$ by \Cref{Goodearl 2.24}.
If $\iota:R \rightarrow P$ and $\mu:R \rightarrow E$ are the inclusion maps, then by injectivity of $E$, there exists a homomorphism $\nu: P \rightarrow E$ such that $\nu\iota=\mu$.
Observe that $R \cap \ker{\nu}=\ker{\mu}=0$.
This implies $\ker{\nu}=0$ since $R$ is essential in $P$ and $\ker{\nu}$ is a submodule of $P$.
Therefore, the inclusion map $\mu:R \rightarrow E$ can be extended to a monomorphism $\nu:P \rightarrow E$.
Moreover, \cite[Theorem 2.26]{Goodearl} shows that $\nu P$ is contained in $Q$, and hence the inclusion map $R \rightarrow Q$ can be extended to a monomorphism $\nu:P \rightarrow Q$.
Finally, note that since $R \leq \nu P \leq Q$ and $R_{R} \leq^{r} Q_{R}$, $Hom_{R}(\nu P \slash R, Q)=0$.
Hence, given $x \in P$, the homomorphism $\sigma:\nu P \slash R \rightarrow Q$ defined by $\sigma(\nu y + R)= \nu(xy)-(\nu x)(\nu y)$ is the zero map.
Therefore, $\nu$ is a ring homomorphism and $Q$ is a maximal right ring of quotients of $R$.
\end{proof}
\end{theorem}

Goodearl shows in \cite[Corollary 2.31]{Goodearl} that $Q^r$ is injective as a right $R$-module. 
Therefore, $Q^r(R)$ is an injective hull of $R$ since $R_{R} \leq^{e} Q^{r}_{R}$ by \Cref{Goodearl 2.24}.
Moreover, since the injective hull is unique up to isomorphism, we can refer to $Q^r(R)$ as the injective hull of $R$.
The following results about maximal quotient rings will be needed later.
The proofs are omitted.

\begin{prop}\cite[Proposition 2.2]{Albrecht}\label{ADF 2.2}
For a right non-singular ring $R$, $R$ is a left p.p.-ring such that $Q^r(R)$ is torsion-free as a right R-module if and only if all non-singular right R-modules are torsion-free.
\end{prop}

\begin{theorem}\cite[Ch. XII, Proposition 7.2]{Stenstrom}\label{NS FG embedding}
If $R$ is a right non-singular ring and $M$ is a finitely generated non-singular right R-module, then there exists a monomorphism 
\\$\varphi: M \rightarrow \oplus_{n} Q^r$ for some $n < \omega$.
In other words, $M$ is isomorphic to a submodule of a free $Q^r$-module.
\end{theorem}

For a ring $R$, its maximal right ring of quotients $Q^{r}$ is a {\it perfect left localization of $R$} if $Q^{r}$ is flat as a right $R$-module and the multiplication map $\varphi: Q^{r} \bigotimes_{R} Q^{r} \rightarrow Q^{r}$, defined by $\varphi(a \otimes b) = ab$, is an isomorphism.  
If $R$ is a right non-singular ring for which $Q^{r}$ is a perfect left localization, then $R$ is called {\it right strongly non-singular}.
Goodearl provides the following useful characterization of right strongly non-singular rings:

\begin{theorem}\cite[Theorem 5.17]{Goodearl}\label{Qr perf rt local}
Let R be a right non-singular ring.  Then, R is right strongly non-singular if and only if every finitely generated non-singular right R-module is isomorphic to a finitely generated submodule of a free right R-module.  
\end{theorem}

\begin{corollary}\cite[Theorem 5.18]{Goodearl}\label{Goodearl 5.18}
Let $R$ be a right non-singular ring.  Then, $R$ is right semi-hereditary, right strongly non-singular if and only if every finitely generated non-singular right $R$-module is projective.
\begin{proof}
For a right non-singular ring $R$, suppose $R$ is right semi-hereditary, right strongly non-singular.
Let $M$ be a finitely generated non-singular right $R$-module.
By \Cref{Qr perf rt local}, $M$ is isomorphic to a finitely generated submodule of a free right $R$-module $F$.
Therefore, since $R$ is right semi-hereditary, $M$ is projective by \Cref{Rot 4.30}.

Conversely, assume every finitely generated non-singular right $R$-module is projective.
Since $R$ is right non-singular, every finitely generated right ideal of $R$ is non-singular.
Hence, every finitely generated right ideal is projective and $R$ is right semi-hereditary.
Furthermore, every finitely generated non-singular right $R$-module is a direct summand, and hence a submodule, of a free right $R$-module.
Therefore, $R$ is right strongly non-singular by \Cref{Qr perf rt local}.
\end{proof}
\end{corollary}

\section{Coincidence of Classes of Torsion-free and Non-singular Modules}

We know turn our attention to rings for which the classes of torsion-free and non-singular right $R$-modules coincide, which is investigated in \cite{Albrecht} by Albrecht, Dauns, and Fuchs. 
A few definitions are needed before stating their theorems in full.
A ring is {\it right semi-simple} if it can be written as a direct sum of modules which have no proper nonzero submodules, and a ring is {\it right Artinian} if it satisfies the descending chain condition on right ideals.
Assume {\it semi-simple Artinain} to mean {\it right semi-simple, right Artinian}.
The following results from Stenstr\"{o}m consider rings with semi-simple right maximal ring of quotients.

\begin{prop}\cite[Ch. XI, Proposition 5.4]{Stenstrom}\label{Qr=Ql iff flat} %XI 5.4
Let R be a ring whose maximal right ring of quotients is semi-simple. Then, $Q^r = Q^l$ if and only if $Q^r$ is flat as a right R-module.
\end{prop}

\begin{theorem}\cite[Ch. XII, Corollaries 2.6,2.8]{Stenstrom}\label{Q ss then left SNS} %XII 2.8
Let $R$ be a ring and suppose $Q^r(R)$ is semi-simple. Then:
\begin{enumerate}[(a)]
\item $Q^r$ is a perfect right localization of R.  In other words, if $R$ is left non-singular, then it is left strongly non-singular.
\item If $M$ is any non-singular right $R$-module, then $M \bigotimes_{R} Q^r$ is the injective hull of $M$.
\end{enumerate}  
\end{theorem}

A ring $R$ is {\it von Neumann regular} if, given any $r \in R$, there exists some $s \in R$ such that $r=rsr$.
These rings are of interest because $R$ is von Neumann regular if and only if every right $R$-module is flat \cite[Theorem 4.9]{Rotman}. 
The following lemmas will be needed in the next chapter. 

\begin{lemma}\cite{Rotman}\label{semi reg}
If R is a semi-simple Artinian ring, then R is von Neumann regular.
\begin{proof}
The Wedderburn-Artin Theorem states that $R$ is semi-simple Artinian if and only if it is isomorphic to a finite direct product of matrix rings over division rings.
For any division ring $D$, $Mat_{n}(D) \cong End_{D}(\bigoplus^n D)$ is von Neumann regular \cite{Rotman}.
Therefore, $R$ is von Neumann regular since direct products of regular rings are regular.
\end{proof}
\end{lemma}

\begin{lemma}\cite{Stenstrom}\label{nonsing reg} %XII 2.2
A ring R is right non-singular if and only if $Q^r$ is von Neumann regular.
\begin{proof}
Stenstr\"{o}m shows in \cite[Ch. XII]{Stenstrom} that if $R$ is right non-singular, then $Q^r \cong End_{R}(E)$, where $E \cong Q^r$ is the injective hull of $R$.
In \cite[Ch. V, Proposition 6.1]{Stenstrom}, it is shown that such rings are regular.

Conversely, assume $Q^r$ is von Neumann regular.  
Let $I$ be an essential right ideal of $R$ and take $x \in R$ to be nonzero.
Suppose $xI=0$.
Since $Q^r$ is regular, there exists some $q \in Q$ such that $xqx=x$.
Hence, $qxR$ is a nonzero right ideal of $R$, and so $I \cap qxR \neq 0$.
Thus, $0 \neq qxr \in I$ for some nonzero $r \in R$.
However, $xr=xqxr \in xI = 0$.
This implies $qxr=0$, which is a contradiction.
Therefore, $xI \neq 0$ and $R$ is right non-singular.
\end{proof}
\end{lemma}

Let $R$ be a ring and $M$ a right $R$-module.
A submodule $U$ of $M$ is {\it $\textbf{S}$-closed} if $M \slash U$ is non-singular. The following lemma shows that annihilators of elements are $\textbf{S}$-closed for non-singular rings.

\begin{lemma}\label{NS ann Sclosed}
If $R$ is a right non-singular ring, then for any $x \in R$, $ann_{r}(x)$ is $\textbf{S}$-closed.
\begin{proof}
Let $R$ be right non-singular.  It needs to be shown that $R/ann_{r}(x)$ is non-singular for any $x \in R$.  
That is, for $x \in R$, $\textbf{Z}(R/ann_{r}(x)) = \{r + ann_{r}(x)$ $|$ $(r + ann_{r}(x))I = 0$ for some $I \leq^{e} R\} = 0$.
Let $0 \neq r + ann_{r}(x) \in R/ann_{r}(x)$ and $I$ be a nonzero essential right ideal of $R$ such that $(r +ann_{r}(x))I = 0$.
Then, for any $a \in I$, $ra + ann_{r}(x) = 0$.
Hence, $ra \in ann_{r}(x)$ and $xra = 0$ for every $a \in I$.
In other words, $(xr)I = 0$.
If $xr \neq 0$, then there is a contradiction since $I \leq^{e} R$ and $\textbf{Z}(R) = 0$.
Thus, $xr = 0$ and $r \in ann_{r}(x)$.
Therefore, $r + ann_{r}(x) = 0$, and it follows readily that $\textbf{Z}(R/ann_{r}(x)) = 0$.
\end{proof}
\end{lemma}

If $R$ is a right non-singular ring and every $\textbf{S}$-closed right ideal of $R$ is a right annihilator, then $R$ is referred to as a {\it right Utumi ring}.
Similarly, $R$ is a {\it left Utumi ring} if $R$ is left non-singular and every $\textbf{S}$-closed left ideal of $R$ is a left annihilator.
The following result from Goodearl characterizes non-singular rings which are both right and left Utumi.

\begin{theorem}\cite[Theorem 2.38]{Goodearl}\label{Ql = Qr} %2.38
If R is a right and left non-singular ring, then $Q^r = Q^l$ if and only if every $R$ is both right and left Utumi.
\end{theorem}

For a ring $R$, if every direct sum of nonzero right ideals of $R$ contains only finitely many direct summands, then $R$ is said to have {\it finite right Goldie-dimension}.
Denote the Goldie-dimension of $R$ as $G$-dim $R_{R}$.
If a ring $R$ with finite right Goldie-dimension also satisfies the ascending chain condition on right annihilators, then $R$ is a {\it right Goldie-ring}.
The maximal right quotient ring $Q^r$ is a {\it semi-perfect left localization of $R$} if $Q^{r}_{R}$ is torsion-free and the multiplication map $Q^r \bigotimes_{R} Q^{r} \rightarrow Q^r$ is an isomorphism.
The following is a useful characterization of rings with finite right Goldie dimension:

\begin{theorem}\cite[Ch. XII, Theorem 2.5]{Stenstrom}\label{Q semi Goldie dim} %CH. XII 2.5
If $R$ is a right non-singular ring, then $Q^r$ is semi-simple if and only if R has finite right Goldie dimension.
\end{theorem}

We are now ready to state two key results form Albrecht, Fuchs, and Dauns, which consider rings for which the classes of torsion-free and non-singular modules coincide.  
These will be needed in the next chapter to prove the main theorem of this thesis.
The proof of \Cref{TF and NS 3.7} is omitted.

\begin{theorem}\cite[Theorem 3.7]{Albrecht}\label{TF and NS 3.7}
The following are equivalent for a ring R:
\begin{enumerate}[(a)]
   \item R is a right Goldie right p.p.-ring and $Q^r$ is a semi-perfect left localization of R.
   \item R is a right Utumi p.p.-ring which does not contain an infinite set of orthogonal idempotents.
   \item R is a right non-singular ring which does not contain an infinite set of orthogonal idempotents, and every finitely generated non-singular right R-module is torsion-free. 
   \item A right R-module M is torsion-free if and only if M is non-singular.
\end{enumerate}
  
\noindent Furthermore, if R satisfies any of the equivalent conditions, then R is a Baer-ring and $Q^r$ is semi-simple Artinian.
\end{theorem}

\begin{theorem}\cite{Albrecht}\label{ADF 4.2}
The following are equivalent for a ring R:
\begin{enumerate}[(a)]
   \item R is a right and left non-singular ring which does not contain an infinite set of orthogonal idempotents, and every $\textbf{S}$-closed left or right ideal is generated by an idempotent.
   \item R is a right or left p.p.-ring, and $Q^r = Q^l$ is semi-simple Artinian.  
   \item R is a right strongly non-singular right p.p.-ring which does not contain an infinite set of orthogonal idempotents.
   \item R is right strongly non-singular, and a right R-module is torsion-free if and only if it is non-singular.
   \item For a right R-module M, the following are equivalent:
   \begin{enumerate}[(i)]
      \item M is torsion-free
      \item M is non-singular
      \item If $E(M)$ is the injective hull of $M$, then $E(M)$ is flat.
     \end{enumerate} 
  \end{enumerate}
  
\begin{proof}
($a$) $\Rightarrow$ ($b$): Assume $R$ is right and left non-singular, contains no infinite set of orthogonal idempotents, and every $\textbf{S}$-closed right or left ideal is generated by an idempotent.
Let $I$ be an $\textbf{S}$-closed right ideal of $R$.
Then, $I = eR$ for some idempotent $e \in R$.  
As shown in the proof of \Cref{eR is ann}, $eR = ann_{r}(1-e)$.
Thus, $I = eR$ is the right annihilator of $1-e$. 
Note that a symmetric argument shows that if $J$ is an $\textbf{S}$-closed left ideal of $R$, then $J = Rf$ is a left annihilator of $1-f$ for some idempotent $f \in R$.
Hence, $R$ is both a right and left Utumi ring.
By \Cref{NS ann Sclosed}, since $R$ is a right non-singular ring, $ann_{r}(x)$ is $\textbf{S}$-closed for every $x \in R$.  
This implies that $ann_{r}(x)$ is generated by an idempotent for every $x \in R$.  
Therefore, $R$ is a right p.p-ring.
A symmetric argument shows that $R$ is also a left p.p.-ring since condition ($a$) applies to both right and left ideals.
Note that $R$ satisfies condition ($b$) of \Cref{TF and NS 3.7} since it is a right Utumi p.p.-ring which does not contain an infinite set of orthogonal idempotents.  
Hence, $Q^r$ is semi-simple Artinian by \Cref{TF and NS 3.7}.
Furthermore, since every right and left $\textbf{S}$-closed ideal is an annihilator, $R$ is right and left Utumi.
Therefore, $Q^r = Q^l$ by \Cref{Ql = Qr}.

($b$) $\Rightarrow$ ($c$): Suppose $R$ is a right p.p.-ring and $Q^r = Q^l$ is semi-simple Artinian.
Since $R$ is a right p.p.-ring, it is also a right non-singular ring.
Hence, $R$ has finite right Goldie dimension by \Cref{Q semi Goldie dim}.
Suppose $R$ contains an infinite set of orthogonal idempotents.  
Consider two orthogonal idempotents $e$ and $f$, and let $x \in eR$ $\cap$ $fR$.
Then, $x = er = fs$ for some $r$, $s \in r$.
This implies that $x = 0$ since $er = e^2r = efs = 0$.
Thus, $eR$ $\cap$ $fR = 0$ for any two orthogonal idempotents $e$ and $f$ in the infinite set, and $eR \bigoplus fR$ is direct.
Hence, $R$ contains an infinite direct sum of nonzero right submodules, which contradicts $R$ having finite right Goldie dimension.
Therefore, $R$ does not contain an infinite set of orthogonal idempotents.

By \Cref{Q ss then left SNS}, since $R$ is semi-simple Artinian, $R$ is a left strongly non-singular ring.
Hence, the multiplication map $\varphi: Q^r \bigotimes_{R} Q^r \rightarrow Q^r$, defined by $\varphi(q \otimes p) = qp$, is an isomorphism.
Note that this also implies that $Q^r$ is flat as a left R-module.
However, in order for $R$ to be right strongly non-singular, it needs to be shown that $Q^r$ is flat as a right R-module.
By \Cref{Qr=Ql iff flat}, $Q^r$ is indeed flat as a right R-module since $Q^r = Q^l$ is assumed to be semi-simple Artinian.
Therefore, $R$ is a right strongly non-singular ring which does not contain an infinite set of orthogonal idempotents.
Note that \Cref{RwCC 8.4} shows that $R$ is also a left p.p.-ring. 
Thus, if we had instead assumed that $R$ is a left p.p.-ring, then a symmetric argument could be used to show that $R$ is also a right p.p.-ring, and the latter part of the proof would remain the same.

($c$) $\Rightarrow$ ($d$): Assume $R$ is a right strongly non-singular right p.p.-ring which does not contain an infinite set of orthogonal idempotents.
Then, $Q^r$ is flat as a right R-module, which follows from $R$ being right strongly non-singular.
Since flat R-modules are torsion-free, this implies that $Q^r$ is torsion-free.
By \Cref{RwCC 8.4}, since $R$ is a right p.p.-ring and does not contain an infinite set of orthogonal idempotents, $R$ is also a left p.p.-ring.
Hence, every non-singular right R-module is torsion-free by \Cref{ADF 2.2}.
Thus, $R$ satisfies condition ($c$) of \Cref{TF and NS 3.7}, which implies that a right R-module M is torsion-free if and only if M is non-singular.

($d$) $\Rightarrow$ ($e$): Suppose $R$ is right strongly non-singular, and a right R-module is torsion-free if and only if it is non-singular.
Then, conditions ($i$) and ($ii$) of ($e$) are clearly equivalent, and it suffices to show that a right R-module is non-singular if and only if its injective hull is flat.
Suppose $M$ is a non-singular right R-module.
Note that $R$ satisfies condition ($d$) of \Cref{TF and NS 3.7}, and hence $Q^r$ is semi-simple Artinian.
By \Cref{Q ss then left SNS}, $M \bigotimes_{R} Q^r$ is an injective hull of $M$.
Thus, if $E(M)$ denotes the injective hull of $M$, then $E(M) \cong M \bigotimes_{R} Q^r$, since an injective hull of a right R-module is unique up to isomorphism.
This implies that $E(M)$ is a right $Q^r$-module, since $M \bigotimes_{R} Q^r$ is a right $Q^r$-module by \Cref{Rot 2.51}.
Furthermore, since $Q^r$ is semi-simple Artinian, every $Q^r$-module is projective.
Hence, $E(M)$ is projective and thus isomorphic to a direct summand of a free $Q^r$-module $F$.
Note that $Q^r$ is flat as a right R-module since $R$ is right strongly non-singular.  
Thus, \Cref{direct sum of flat} shows that any free $Q^r$-module is flat since such modules can be written as $\bigoplus_{i \in I} M_{i}$ for some index set $I$, where $M_{i}$ is isomorphic to $Q^r$ for every $i \in I$.
This implies that $E(M)$ is flat by \Cref{direct sum of flat} since it is a direct summand of the flat right $R$-module $F = \bigoplus_{i \in I} M_{i}$.

On the other hand, assume that the injective hull $E(M)$ of some right R-module $M$ is flat.  
Noting again that $R$ satisfies condition ($d$) of \Cref{TF and NS 3.7}, it follows that $R$ is a right p.p.-ring.
Thus, $R$ is a torsion-free ring by \Cref{Dauns Fuchs 4.5}.
Since flat $R$-modules are torsion-free, $E(M)$ is torsion-free as a right R-module. 
Hence, $M$ is a submodule of a torsion-free right R-module.  
Thus, $M$ is a torsion free right R-module by \Cref{submod of tor free ring}.  
Therefore, $M$ is non-singular since a right $R$-module is torsion-free if and only if it is non-singular by assumption.

($e$) $\Rightarrow$ ($a$): For a right $R$-module $M$, assume that $M$ is torsion-free if and only if $M$ is non-singular if and only if the injective hull $E(M)$ is flat.  
By \Cref{TF and NS 3.7}, $R$ is a right p.p.-ring which does not contain an infinite set of orthogonal idempotents.
It then follows from \Cref{pp implies ns} that $R$ is a right non-singular ring.
Hence, $R$ is also a left p.p.-ring by \Cref{ADF 2.2}, since every non-singular right R-module is torsion-free, and a symmetric argument for \Cref{pp implies ns} shows that $R$ is left non-singular.

The injective hull $E(R)$ is flat as a right R-module since $R$ is assumed to be right non-singular.
Hence, $Q^r$ is flat as a right R-module, since $Q^r$ is the injective hull of $R$.
We've already shown that $R$ satisfies the equivalent conditions of \Cref{TF and NS 3.7}, which implies that $Q^r$ is a semi-simple Artinian ring. 
Thus, it follows from \Cref{Qr=Ql iff flat} that $Q^r = Q^l$.
Since $R$ is both right and left non-singular, every $\textbf{S}$-closed right ideal of $R$ is a right annihilator and every $\textbf{S}$-closed left ideal of $R$ is a left annihilator by \Cref{Ql = Qr}.  
Furthermore, note that \Cref{TF and NS 3.7} shows that $R$ is a Baer-ring.
Hence, every annihilator is generated by an idempotent.
Therefore, every $\textbf{S}$-closed right ideal and every $\textbf{S}$-closed left ideal is generated by an idempotent. 
\end{proof}
\end{theorem}

\chapter{Morita Equivalence}
Before proving the main theorem, we discuss Morita equivalences.
In particular, we show that there is a Morita equivalence between $R$ and $Mat_{n}(R)$ for any $0 < n < \omega$.
This is then used to show that the classes of torsion-free and non-singular $Mat_{n}(R)$-modules coincide for certain conditions placed on $R$.

Let R and S be rings.  
The categories $Mod_{R}$ and $Mod_{S}$ are {\it equivalent} (or {\it isomorphic}) if there are functors $F: Mod_{R} \rightarrow Mod_{S}$ and $G: Mod_{S} \rightarrow Mod_{R}$ such that $FG \cong 1_{Mod_{S}}$ and $GF \cong 1_{Mod_{R}}$.
Note that these are natural isomorphisms.
In other words, if $\eta:GF \rightarrow 1_{Mod_{R}}$ denotes the natural isomorphism, then for each $M,N \in Mod_{R}$, there exist isomorphisms $\eta_{M}:GF(M) \rightarrow M$ and $\eta_{N}:GF(N) \rightarrow N$ such that $\beta\eta_{M}=\eta_{N}GF(\beta)$ whenever $\beta \in Hom_{R}(M,N)$.
Here, $GF(\beta)$ denotes the induced homomorphism.
The functors $F$ and $G$ are referred to as an {\it equivalence} of $Mod_{R}$ and $Mod_{S}$.
If such an equivalence exists, then R and S are said to be {\it Morita-equivalent}. 
In \cite[Ch. IV, Corollary 10.2]{Stenstrom}, Stenstr\"{o}m shows that R and S are Morita-equivalent if and only if there are bimodules $_{S}P_{R}$ and $_{R}Q_{S}$ such that $P \bigotimes_{R} Q \cong S$ and $Q \bigotimes_{S} P \cong R$.
A property $P$ is referred to as {\it Morita-invariant} if for every ring $R$ satisfying $P$, every ring $S$ Morita-equivalent to $R$ also satisfies $P$.

A {\it generator} of $Mod_{R}$ is a right $R$-module $P$ satisfying the condition that every right $R$-module $M$ is a quotient of $\displaystyle{ \bigoplus}_{I} \hspace{.1cm}P\hspace{.1cm}$. Note that $R$ and any free right $R$-module are generators of $Mod_{R}$.
A {\it progenerator} of $Mod_{R}$ is a generator which is finitely generated and projective.

\begin{lemma}\cite{Anderson}\label{Equivalence}
Let $R$ be a ring, $P$ a progenerator of $Mod_{R}$, and $S = End_{R}(P)$.  Then, there is an equivalence $F: Mod_{R} \rightarrow Mod_{S}$ given by $F(M) = Hom_{R}(P,M)$ with inverse $G: Mod_{S} \rightarrow Mod_{R}$ given by $G(N) = N \bigotimes_{S} P$.
\begin{proof}
As a projective generator of $Mod_{R}$, $P$ is a right $R$-module.
$P$ also has a left $S$-module structure with $(f*g)(x)=f(g(x))$ for $f, g \in S$ and $x \in P$, where multiplication in the endomorphism ring is defined as composition of functions.
It then readily follows that $P$ is an $(S,R)$-bimodule since $f(xr)=f(x)r$ for any $f \in S$ and $r \in R$.
Thus, $F = Hom_{S}(P,\underline{\hspace{.3cm}})$ is a functor $Mod_{S} \rightarrow Mod_{R}$ and $G = \underline{\hspace{.3cm}} \bigotimes_{R} P$ is a functor $Mod_{R} \rightarrow Mod_{S}$ by \Cref{Functor into ModR ModS}.

It needs to be shown that $GF \cong 1_{Mod_{R}}$ and $FG \cong 1_{Mod_{S}}$ are natural isomorphisms.
Since $P$ is a progenerator of $Mod_{R}$, it is finitely generated and projective as a right $R$-module.
Thus, it follows from \Cref{And 20.11} that if $M$ is any right $R$-module, then $GF(M)=G(Hom_{R}(P,M))=Hom_{R}(P,M) \bigotimes_{S} P \cong Hom_{R}(Hom_{S}(P,P),M) \cong Hom_{R}(End_{S}(P),M) \\ \cong Hom_{R}(R,M) \cong M$.
Similarly, given any right $S$-module $N$, $FG(N)=F(N \bigotimes_{S} P)=Hom_{R}(P,N \bigotimes_{S} P) \cong N \bigotimes_{S} Hom_{R}(P,P) = N \bigotimes_{S} S \cong N$ by \Cref{And 20.10}.
Therefore, $F$ is an equivalence with inverse $G$.
\end{proof}
\end{lemma}

\begin{prop}\label{mat morita}
Let $R$ be a ring.  For every $0 < n < \omega$, $R$ is Morita-equivalent to $Mat_{n}(R)$.
\begin{proof}
Let $P$ be a finitely generated free right $R$-module with basis $\{x_{i}\}_{i=1}^{n}$ for $0 < n < \omega$.
Then, $P$ is a progenerator of $Mod_{R}$ and $Mat_{n}(R) \cong End_{R}(P)$ by \Cref{End iso Mat}. 
Therefore, the equivalence of \Cref{Equivalence} is a Morita-equivalence between $R$ and $Mat_{n}(R)$.
\end{proof}
\end{prop}

\begin{lemma}\cite[Ch. X, Proposition 3.2]{Stenstrom}\label{Qr Morita}
If $R$ and $S$ are Morita-equivalent, then the maximal ring of quotients, $Q^r(R)$ and $Q^r(S)$, are also Morita equivalent.
\end{lemma}

%\begin{theorem}\cite{Albrecht}\label{ADF 5.1}
%The following are equivalent for a ring R:
%\begin{enumerate}[(a)]
 %  \item R is right strongly non-singular, right semi-hereditary, and does not contain an infinite set of orthogonal idempotents.
 %  \item If S is Morita equivalent to R, then S is right strongly non-singular, and the classes of torsion-free right $S$-modules and non-singular right $S$-modules coincide.   
 % \end{enumerate}
%\begin{proof}
%($a$) $\Rightarrow$ ($b$): 
%
%($b$) $\Rightarrow$ ($a$): Assume every ring $S$ which is Morita-equivalent to $R$ is right strongly non-singular, and a right S-module is torsion-free if and only if it is non-singular.
%Since $R$ is Morita-equivalent with itself, $R$ is right strongly non-singular and an $R$-module is torsion-free if and only if it is non-singular.
%Thus, $R$ satisfies condition ($d$) of \Cref{ADF 4.2}, which implies that $R$ does not contain an infinite set of orthogonal idempotents.
%To see that $R$ is right semi-hereditary, consider $Mat_{n}(R)$, which is Morita-equivalent to $R$ for any $0 < n < \omega$ by \Cref{mat morita}.
%By assumption, this implies that every right $Mat_{n}(R)$-module is torsion-free if and only if it is non-singular.
%Hence, $Mat_{n}(R)$ is a right p.p.-ring for every $n$ by \Cref{TF and NS 3.7}.  
%Therefore, $R$ is right semi-hereditary by \Cref{semihered iff Mat pp}.
%Thus, $R$ is a right strongly non-singular, right semi-hereditary ring which does not contain any infinite set of orthogonal idempotents.
%\end{proof}
%\end{theorem}
%
\begin{prop}\label{MI Properties}
Let $R$ and $S$ be Morita-equivalent rings with equivalence 
\\$F: Mod_{R} \rightarrow Mod_{S}$ and $G: Mod_{S} \rightarrow Mod_{R}$.
\begin{enumerate}[(i)]
\item If $U$ is an essential submodule of a right $R$-module $M$, then $F(U)$ is an essential submodule of the right $S$-module $F(M)$.
\item If $M$ is a non-singular right $R$-module, then $F(M)$ is a non-singular right $S$-module.
\end{enumerate}
In other words, essentiality and non-singularity are Morita-invariant properties.
\begin{proof}
($i$): Let $U \leq^{e} M$.
Then, the inclusion map $\iota:U \rightarrow M$ is an essential monomorphism.
Consider the induced homomorphism $F(\iota):F(U) \rightarrow F(M)$.
Note that since $\iota$ is a monomorphism, $F(\iota)$ is a monomorphism \cite[Proposition 21.2]{Anderson}.
Let $W$ be any right $S$-module and take $\beta \in Hom_{S}(F(M),W)$ to be such that $\beta F(\iota):F(U) \rightarrow W$ is a monomorphism.
There is a natural isomorphism $\Phi_{U,W}:Hom_{S}(F(U),W) \rightarrow Hom_{R}(U,G(W))$ defined by $\gamma \mapsto G(\gamma)\eta^{-1}_{U}$, where $\eta_{U}$ denotes the isomorphism $GF(U) \rightarrow U$ \cite[21.1]{Anderson}.
Hence, $\Phi_{U,W}(\beta F(\iota))$ is a monomorphism.
Moreover, $\Phi_{U,W}(\beta F(\iota))= G(hF(\iota))\eta^{-1}_{U} =G(h)GF(\iota)\eta^{-1}_{U}=G(h)\eta_{M}^{-1}\eta_{M} GF(\iota)\eta_{U}^{-1}=\Phi_{M,W}(\beta)\iota \eta_{U}\eta_{U}^{-1}=\Phi_{M,W}(\beta)\iota$.
Thus, $\Phi_{M,W}(\beta)\iota$ is a monomorphism and it follows from \Cref{And 5.13} that $\Phi_{M,W}(\beta)$ is a monomorphism since $\iota$ is essential.
Furthermore, $\Phi_{M,W}(\beta)$ is a monomorphism if and only if $\beta$ is a monomorphism \cite[Lemma 21.3]{Anderson}.
Hence, $F(\iota)$ is an essential monomoprhism by \Cref{And 5.13}. 
Therefore, $F(U) \cong im(F(\iota)) \leq^{e} F(M)$.

($ii$): Let $M$ be a non-singular right $R$-module.
It needs to be shown that $F(M)$ is a non-singular right $S$-module and in view of \Cref{Goodearl 1.20a} it suffices to show that $Hom_{S}(C,F(M))= 0$ for any singular right $S$-module $C$.
By \Cref{Goodearl 1.20b}, there is an exact sequence $0 \rightarrow A \xrightarrow{f} F \rightarrow C \rightarrow 0$ of right $S$-modules such that $f(A) \leq^{e} F$ and $F$ is free.
Then, $G(f(A)) \leq^{e} G(B)$ by ($i$).
Hence, $0 \rightarrow G(A) \xrightarrow{G(f)} G(B) \rightarrow G(C) \rightarrow 0$ is an exact sequence of right $R$-modules such that $G(f(A)) \leq^{e} G(B)$.
Thus, $G(C)$ is a singular right $R$-module by \Cref{Goodearl 1.20b}.
Since $G(C)$ is singular and $M$ is non-singular, $Hom_{R}(G(C),M) = 0$ by \Cref{Goodearl 1.20a}.
Therefore, $Hom_{S}(C,F(M)) \cong Hom_{R}(G(C),M) = 0$.
Observe that in this proof, it is also shown that singularity is Morita-invariant since we show that $G(C)$ is singular for an arbitrary singular module $C$.
\end{proof}
\end{prop}

We now prove the main theorem of this thesis.

\begin{theorem}\label{main}
The following are equivalent for a ring $R$:
\begin{enumerate}[(a)]
\item $R$ is a right strongly non-singular, right semi-hereditary, right Utumi ring not containing an infinite set of orthogonal idempotents.
\item Whenever $S$ is Morita-equivalent to $R$, then the classes of torsion-free right $S$-modules and non-singular right $S$-modules coincide.
\item For every $0 < n < \omega$, $Mat_{n}(R)$ is a right and left Utumi Baer-ring not containing an infinite set of orthogonal idempotents.
\end{enumerate}
Moreover, if $R$ is such a ring, then the corresponding left conditions are also satisfied.
\begin{proof}
($a$) $\Rightarrow$ ($b$): Assume R is a right strongly non-singular, right semi-hereditary, right Utumi ring not containing an infinite set of orthogonal idempotents.
Let $R$ and $S$ be Morita equivalent, and let $F: Mod_{R} \rightarrow Mod_{S}$ and $G: Mod_{S} \rightarrow Mod_{R}$ be an equivalence.
Also, take $N$ to be a finitely generated non-singular right R-module.
Since $R$ is right strongly non-singular, $N$ is isomorphic to finitely generated submodule $V$ of a free right R-module by \Cref{Qr perf rt local}.
Furthermore, since $R$ is right semi-hereditary and free R-modules are projective, $V \cong N$ is projective by \Cref{Rot 4.30}.
Thus, since projective modules are torsion-free, it follows that finitely generated non-singular right R-modules are torsion-free.
Therefore, R satisfies condition ($c$) of \Cref{TF and NS 3.7}, which implies that the maximal ring of quotients $Q^r(R)$ is a semi-simple Artinian ring.
Note that $Q^r(R)$ and $Q^r(S)$ are Morita-equivalent by \Cref{Qr Morita}.
Hence, $Q^r(S)$ is also semi-simple Artinian, since semi-simpleness and Artinian are properties preserved under a Morita-equivalence \cite{Anderson}.
Furthermore, $Q^r(S)$ is a regular ring by \Cref{semi reg}.
Therefore, \Cref{nonsing reg} shows that $S$ is right non-singular.

Let $M$ be a finitely generated non-singular right S-module.
Then, $G(M)$ is a finitely generated non-singular right R-module since non-singularity and being finitely generated are both Morita-invariant properties \cite{Anderson}.
Thus, since $R$ is a right strongly non-singular ring, $G(M)$ is isomorphic to a finitely generated submodule of a free right R-module $P$ by \Cref{Qr perf rt local}.
Note that as a free right R-module, $P$ is projective, which is also a Morita-invariant property \cite{Anderson}.
Hence, $F(P)$ is a projective right S-module.
Furthermore, since $G(M)$ is isomorphic to a finitely generated submodule of $P$, $FG(M) \cong M$ is isomorphic to a finitely generated submodule $U$ of $F(P)$.
Now, $F(P)$ is projective and hence a submodule of a free right S-module, which implies $U \cong M$ is a submodule of a free right $S$-module.
Therefore, $M$ is isomorphic to a finitely generated submodule of a free right S-module, and $S$ is right strongly non-singular by \Cref{Qr perf rt local}.

It has been shown that $S$ is a right non-singular ring with a semi-simple Artinian maximal right ring of quotients.
Thus, $S$ has finite right Goldie dimension by \Cref{Q semi Goldie dim}.
Hence, $S$ cannot contain an infinite set of orthogonal idempotents.  
Moreover, $S$ is a right p.p.-ring since $R$ is right semi-hereditary.
For if $P$ is a principal right ideal of $S$, then $G(P)$ is a finitely generated right ideal of the right semi-hereditary ring $R$, which implies that $G(P)$ is projective.
Hence, $FG(P) \cong P$ is projective, which again follows from projectivity being Morita-invariant.
Then, $S$ is a right strongly non-singular right p.p.-ring which does not contain an infinite set of orthogonal idempotents.
Therefore, a right S-module is torsion-free if and only if it is non-singular by \Cref{ADF 4.2}.

($b$) $\Rightarrow$ ($a$): Assume that the classes of torsion-free and non-singular $S$-modules coincide for every ring $S$ Morita-equivalent to $R$.
Thus, since $Mat_{n}(R)$ is Morita-equivalent to $R$ for every $0<n<\omega$, the classes of torsion-free right $Mat_{n}(R)$-modules and non-singular right $Mat_{n}(R)$-modules coincide for every $0<n<\omega$.
Hence, $Mat_{n}(R)$ is a right Utumi p.p.-ring which does not contain an infinite set of orthogonal idempotents by \Cref{TF and NS 3.7}.
Thus, $R$ is right semi-hereditary by \Cref{semihered iff Mat pp}. 
In particular, since these conditions are satisfied for every $0<n<\omega$, they are satisfied for $n=1$. 
Hence, $R \cong Mat_{1}(R)$ is a right semi-hereditary right Utumi ring not containing an infinite set of orthogonal idempotents.

It needs to be shown that $R$ is right strongly non-singular.
Let $M$ be a finitely generated non-singular right $R$-module.
By \Cref{Goodearl 5.18}, $R$ is right strongly non-singular if $M$ is projective.
Let $0 \rightarrow U \rightarrow F = \displaystyle \bigoplus^{n} R \rightarrow M \rightarrow 0$ be an exact sequence of right $R$-modules.
Since $F$ is a finitely generated free right $R$-module, it is a progenerator of $Mod_{R}$.  
Hence, $0 \rightarrow Hom_{R}(F,U) \rightarrow Hom_{R}(F,F) = End_{R}(F) \rightarrow Hom_{R}(F,M) \rightarrow 0$ is exact by \Cref{Rotman 3.2}. 
Moreover, if $S = End_{R}(F) \cong Mat_{n}(R)$, then $F: Mod_{R} \rightarrow Mod_{S}$ given by $F(M) = Hom_{R}(F,M)$ and $G: Mod_{S} \rightarrow Mod_{R}$ given by $G(N) = N \bigotimes_{S} F$ is an equivalence by \Cref{Equivalence}.
Thus, $Hom_{R}(F,M)$ is a non-singular right $S$-module by \Cref{MI Properties} ($ii$). % since $M$ is a non-singular right $R$-module.
Furthermore, since $S$ is Morita-equivalent to $R$, the $S$-module $Hom_{R}(F,M)$ is torsion-free by assumption.
Note that since the sequence is exact, $Hom_{R}(F,M) \cong S \slash Hom_{R}(F,U)$.
Thus, $Hom_{R}(F,M)$ is cyclic as an $S$-module since $Hom_{R}(F,U)$ is a right ideal of the right $S$-module $S$.
Note also that $S$ is a left p.p.-ring by \Cref{RwCC 8.4} since $S$ is a right p.p.-ring which does not contain an infinite set of orthogonal idempotents.
Thus, the cyclic torsion-free right $S$-module $Hom_{R}(F,M)$ is projective by \Cref{AFD 3.4}.
Therefore, $M \cong GF(M) = G(Hom_{R}(F,M))$ is a projective right $R$-module and $R$ is right strongly non-singular.

($a$) $\Rightarrow$ ($c$): Assume $R$ is right strongly non-singular, right semi-hereditary, right Utumi, and does not contain an infinite set of orthogonal idempotents.
It has been shown that any ring $S$ Morita-equivalent to such a ring is right strongly non-singular and the classes of torsion-free and non-singular right $S$-modules coincide.
Thus, $Mat_{n}(R)$ is right strongly non-singular and a right $Mat_{n}(R)$-module is torsion-free if and only if it is non-singular, which follows from $Mat_{n}(R)$ being Morita-equivalent to $R$ for any $0 < n < \omega$.
By \Cref{ADF 4.2}, $Mat_{n}(R)$ is a right strongly non-singular right p.p.-ring which does not contain an infinite set of orthogonal idempotents.
It then follows from \Cref{RwCC 8.4} that $Mat_{n}(R)$ satisfies the ascending chain condition on right annihilators.
Furthermore, \Cref{Dauns Fuchs 4.5} shows that since $Mat_{n}(R)$ is a right p.p.-ring, $Mat_{n}(R)$ is a torsion-free ring such that right annihilators of elements are finitely generated. 
Hence, $Mat_{n}(R)$ is a Baer-ring by \Cref{Fuchs and Dauns 4.8}.
Moreover, \Cref{ADF 4.2} shows that every $\textbf{S}$-closed one-sided ideal of $Mat_{n}(R)$ is generated by an idempotent.
Thus, every right ideal of $Mat_{n}(R)$ is a right annihilator and every left ideal of $Mat_{n}(R)$ is a left annihilator.
Hence, $Mat_{n}(R)$ is a right and left Utumi ring.

($c$) $\Rightarrow$ ($a$): Suppose $Mat_{n}(R)$ is a right and left Utumi Baer-ring for every $0<n<\omega$ and does not contain an infinite set of orthogonal idempotents.
Then, $Mat_{n}(R)$ is a right p.p.-ring, and so $R$ is right semi-hereditary by \Cref{semihered iff Mat pp}.
Furthermore, since $Mat_{n}(R)$ satisfies these conditions for every $0<n<\omega$, $R \cong Mat_{1}(R)$ is a right and left Utumi Baer-ring not containing an infinite set of orthogonal idempotents.
Thus, every $\textbf{S}$-closed one-sided ideal of $R$ is an annihilator and hence generated by an idempotent.
Therefore, since $R$ is a right and left p.p.-ring and hence right and left non-singular, $R$ is right strongly non-singular by \Cref{ADF 4.2}.
\end{proof}
\end{theorem}

\begin{corollary}
The following are equivalent for a ring R which does not contain an infinite set of orthogonal idempotents:
\begin{enumerate}[(a)]
\item R is a right and left Utumi, right semi-hereditary ring.
\item For every $0 < n < \omega$, $Mat_{n}(R)$ is a Baer-ring, and $Q^r(R)$ is torsion-free as a right R-module.
\end{enumerate}
\begin{proof}
($a$) $\Rightarrow$ ($b$): Suppose $R$ is right and left Utumi and right semi-hereditary.
Then, $R$ is a right p.p.-ring and hence right non-singular.
Moreover, $R$ is a left p.p.-ring by \Cref{RwCC 8.4}, which implies that $R$ is also a left non-singular ring.
Since $R$ is both right and left Utumi, $Q^r(R) = Q^l(R)$ by \Cref{Ql = Qr}.
Furthermore, since $R$ is a right Utumi right p.p.-ring which does not contain an infinite set of orthogonal idempotents, $Q^r(R) = Q^l(R)$ is semi-simple Artinian and torsion-free by \Cref{TF and NS 3.7}.
Therefore, $R$ is right strongly non-singular by \Cref{ADF 4.2}. 

Since $R$ is a right strongly non-singular, right semi-hereditary, right Utumi ring not containing an infinite set of orthogonal idempotents, the classes of torsion-free and non-singular right $Mat_{n}(R)$-modules coincide by \Cref{main}.
Moreover, the proof of \Cref{main} shows that $Mat_{n}(R)$ is right strongly non-singular.
Thus, $Mat_{n}(R)$ is a right strongly non-singular, right p.p.-ring not contain an infinite set of orthogonal idempotents by \Cref{ADF 4.2}.
It then follows from \Cref{RwCC 8.4} that $Mat_{n}(R)$ satisfies the ascending chain condition on right annihilators.
Since $Mat_{n}(R)$ is a right p.p.-ring, \Cref{Dauns Fuchs 4.5} shows that $Mat_{n}(R)$ is a torsion-free ring such that right annihilators of elements are finitely generated. 
Hence, $Mat_{n}(R)$ is a Baer-ring by \Cref{Fuchs and Dauns 4.8}.

($b$) $\Rightarrow$ ($a$): Assume $Mat_{n}(R)$ is a Baer-ring for every $0 < n < \omega$, and $Q^r(R)$ is torsion-free as a right $R$-module.
Since $Mat_{n}(R)$ is a Baer-ring, it is both a right and left p.p.-ring. 
Hence, $R$ is both right and left semi-hereditary by \Cref{semihered iff Mat pp}.  
It then readily follows that $R$ is right and left non-singular.
Note also that $R \cong Mat_{1}(R)$ is a Baer-ring since $Mat_{n}(R)$ is Baer for every $0 < n < \omega$. 
Let $I$ be a proper $\textbf{S}$-closed right ideal of $R$.  
Then, $R/I$ is non-singular as a right R-module.
Furthermore, $R/I$ is cyclic and thus finitely generated.
Hence, $R/I$ is isomorphic to a submodule of a free $Q^r$-module by \Cref{NS FG embedding}. 
Since $Q^r$ is assumed to be torsion-free as a right R-module, it follows from \Cref{ADF 3.3} that $I$ is generated by an idempotent $e \in R$.
Hence, $I = ann_{r}(1-e)$ by \Cref{eR is ann} and $R$ is right Utumi.
Observe that the argument works for $\textbf{S}$-closed left ideals as well, and so $R$ is also left Utumi.
\end{proof}
\end{corollary}

We conclude by considering two examples, the first of which illustrates why the condition of being right semi-hereditary is necessary in the main theorem.
Let $R=\mathbb{Z}[x]$.
As an integral domain, $R$ is a strongly non-singular p.p.-ring not containing an infinite set of orthogonal idempotents \cite[Corollary 3.10]{Albrecht}.
By \Cref{ADF 4.2}, the classes of torsion-free and non-singular right $R$-modules coincide, and by \Cref{TF and NS 3.7} $R$ is right Utumi.
However, $R$ is not semi-hereditary since the ideal $x\mathbb{Z}[x] + 2\mathbb{Z}[x]$ of $\mathbb{Z}[x]$ is not projective.
As seen in the proof of \Cref{semihered iff Mat pp}, this implies $S=Mat_{2}(R)$ is not a right or left p.p.-ring, and hence not a Baer ring.
Therefore, the main theorem does not hold if $R$ is not assumed to be right semi-hereditary.
Moreover, this example shows that the classes of torsion-free and non-singular $S$-modules do not necessarily coincide, even if this holds for $R$ and there is a Morita-equivalence between $R$ and $S$.

Finally, we consider an example from \cite{Chatters} which details a ring with finite right Goldie-dimension but infinite left Goldie-dimension.
In the context of this thesis, this example provides a right Utumi Baer-ring which is not left Utumi.
Let $K=F(y)$ for some field $F$ and consider the endomorphism $f$ of $K$ determined by $y \mapsto y^2$.
The ring we consider is $R=K[x]$ with coefficients written on the right and multiplication defined according to $kx=xf(k)$ for $k \in K$.
Observe that $yx=xy^2$.
It can be shown that $Rx \cap Rxy = 0$, and hence $Rxy \bigoplus Rxyx \bigoplus Rxyx^2 \bigoplus ... \bigoplus Rxyx^{k} \bigoplus ...$ is an infinite direct sum of left ideals of $R$.
Thus, $R$ has infinite left Goldie-dimension.
On the other hand, every right ideal of $R$ is a principal ideal \cite{Chatters}, and thus $R$ is right Noetherian.
Hence, $R$ is a right Goldie-ring.
It then follows from \Cref{TF and NS 3.7} that $R$ is a right Utumi Baer ring and $Q^r$ is semi-simple Artinian.
However, $R$ having infinite left Goldie-dimension but finite right Goldie-dimension implies that $Q^r \neq Q^l$ \cite[Proposition 4.1]{Albrecht}.
Therefore, \Cref{Ql = Qr} shows that $R$ cannot be left Utumi.

\end{document}